\documentclass[10pt,british]{amsart}

\newif\ifpdf
\ifx\pdfoutput\undefined\pdffalse\else\pdfoutput=1\pdftrue\fi\newif\pdf
\ifpdf\relax\else
\usepackage[T1]{fontenc}
\fi

\usepackage{babel,eucal,url,amssymb,stmaryrd,
enumerate,amscd,
}

\usepackage[pagebackref]{hyperref}

\theoremstyle{plain}
\newtheorem{thm}{Theorem}[section]
\newtheorem{lem}[thm]{Lemma}
\newtheorem{pro}[thm]{Proposition}
\newtheorem{co}[thm]{Corollary}

\theoremstyle{definition}
\newtheorem{defn}[thm]{Definition}

\theoremstyle{remark}
\newtheorem{rem}[thm]{Remark}
\newtheorem{example}[thm]{Example}
\newtheorem*{ack}{Acknowledgements}

\newcommand{\Gtwo}{\ifmmode{{\rm G}_2}\else{${\rm G}_2$}\fi}

\newcommand{\LC}{{\nabla^g}}

\date{\today}
\begin{document}

\title[ParaHermitian and ParaQuaternionic manifolds]%
{ParaHermitian and ParaQuaternionic manifolds}

\author{Stefan Ivanov}
\address[Ivanov]{University of Sofia "St. Kl. Ohridski"\\
Faculty of Mathematics and Informatics\\ Blvd. James Bourchier 5\\
1164 Sofia, Bulgaria} \email{ivanovsp@fmi.uni-sofia.bg}

\author{Simeon Zamkovoy}
\address[Zamkovoy]{University of Sofia "St. Kl. Ohridski"\\
Faculty of Mathematics and Informatics\\ Blvd. James Bourchier 5\\
1164 Sofia, Bulgaria} \email{}

\begin{abstract}
A set of canonical parahermitian connections on an almost
paraHermitian manifold is defined. ParaHermitian version of the
Apostolov-Gauduchon generalization of the Goldberg-Sachs theorem
in General Relativity is given. It is proved that the Nijenhuis
tensor of a  Nearly paraK\"ahler manifolds is parallel with
respect to the canonical connection. Salamon's twistor
construction on quaternionic manifold is adapted to the
paraquaternionic case. A hyper-paracomplex structure is
constructed on Kodaira-Thurston (properly elliptic) surfaces as
well as on the Inoe surfaces modeled on $Sol^4_1$. A locally
conformally flat hyper-paraK\"ahler (hypersymplectic) structure
with parallel Lee form on Kodaira-Thurston surfaces is obtained.
Anti-self-dual non-Weyl flat neutral metric on Inoe surfaces
modeled on $Sol^4_1$ is presented. An example of anti-self-dual
neutral metric which is
 not locally conformally
hyper-paraK\"ahler  is constructed.


\end{abstract}

\maketitle
\setcounter{tocdepth}{2}
\tableofcontents

\section{Introduction}
We study the geometry of structures on a differentiable manifold
related to the algebra of paracomplex numbers as well as to the
algebra of paraquaternions together with a naturally associated
metric which is necessarily of neutral signature. These structure
lead to the notion of almost paraHermitian manifold, in even
dimension, as well as to the notion of almost paraquaternionic and
hyper-paraHermitian manifolds in dimensions divisible by four.
Some of these spaces, hyper-paracompex and hyper-paraHermitian
manifolds, become attractive in theoretical physics since they
play a role in string theory \cite{OV,hul,JR,Bar,Hull} and
integrable systems \cite{D1}.

Almost paraHermitian geometry is a topic with many analogies with
the almost Hermitian geometry and also with differences. In the
present note we show that a lot of local and some of the global
results in almost Hermitian manifolds carry over, in the
appropriately defined form, to the case of almost paraHermitian
spaces.

We define a set of canonical paraHermitian connections on an almost paraHermitian
manifold and use them to describe properties of 4-dimensional paraHermitian
and 6-dimensional Nearly paraK\"ahler  spaces.

We present a paraHermitian analogue of the Apostolov-Gauduchon
generalization \cite{AG1} of the Goldberg-Sachs theorem in General
Relativity (see e.g. \cite{PR}) which relates the Einstein
condition to the structure of the positive Weyl tensor in
dimension 4. Namely, we prove
\begin{thm}\label{th4}
Let $(M,g,P)$ be a 4-dimensional paraHermitian manifold. Let $W^+$
be the self-dual part of the Weyl tensor and $\theta$ be the Lee
1-form. The following conditions are equivalent:
\begin{enumerate}
\item[a)] The 2-form $d\theta$ is anti-self-dual, $d\theta^+=0$;
\item[b)] $W_2^+=0$, equivalently, the fundamental 2-form
is an eigen-form of $W^+$;
\item[c)] $(\delta W^+)^- =0$, equivalently,
$(\delta W)(X^{1,0};Y^{1,0},Z^{1,0}) =0$.
\end{enumerate}
\end{thm}
\begin{co}
Assume that the Ricci tensor $\rho$ of a paraHermitian 4-manifold
is $P$-anti-invariant, $\rho(PX,PY)=-\rho (X,Y)$. Then  $d\theta $
is anti-self-dual 2-form, $d\theta^+=0$.

In particular, on a paraHermitian Einstein 4-manifold the
fundamental 2-form is an eigen-form of the positive Weyl tensor.
\end{co}
It turns out that any conformal class of neutral metrics on an
oriented 4-manifold is equivalent to the existence of a local
almost hyper-paracomplex structure, i.e. a collection of
anti-commuting almost complex structure and almost para-complex
structure. Using the properties of the Bismut connection, we
derive that the integrability of the almost hyper-paracomplex
structure leads to the anti-self-duality of the corresponding
conformal class of neutral metrics (Theorem~\ref{inq}). Applying
this result to invariant hyper-paracomplex structure on
4-dimensional Lie groups \cite{ASal,BVuk0} we find explicit
anti-self-dual non Weyl flat neutral metrics on some compact
4-manifolds. Some of these metrics seem to be new.

We apply our considerations to Kodaira-Thurston complex surfaces
modeled on $S^1\times \widetilde{SL(2,\mathbb R)}$ (properly
elliptic surfaces) as well as to the Inoe surfaces modeled on
$Sol^4_1$ in the sense of \cite{Wall}. These surfaces do not admit
any (para) K\"ahler structure \cite{Wall,Biq,Pet}. It is also
known that these surfaces do not support a hyper-complex structure
\cite{Ka,Bo}.

In  contrast, we obtain
\begin{thm}\label{kod-t}
The Kodaira-Thurston surfaces $M=S^1\times\widetilde{(SL(2,\mathbb
R)}/\Gamma)$  admit a hyper-paracomplex structure. The
corresponding hyper-paraHermitian structure  has
$\nabla^g$-parallel Lee form and is locally (not globally)
conformally equivalent to a flat hyper-paraK\"ahler
(hypersympectic) structure.
\end{thm}
\begin{thm}\label{inoe}
The Inoe surfaces modeled on $Sol^4_1$ admit a hyper-paracomplex
structure. The corresponding neutral metric is anti-self-dual
non-Weyl flat. The para-hermitian structure is locally (but not
globally) conformally hyper-paraK\"ahler (hypersymplectic).
\end{thm}

The Inoe surfaces modeled on $Sol^4_1$ are compact solvmanifolds.
A compact 4-dimensional solvmanifold $S$ can be written, up to
double covering, as $G/\Gamma$ where $G$ is a simply connected
solvable Lie group and $\Gamma$ is a lattice of $G$ and all
compact four-dimensional solvmanifolds admitting a complex
structure are classified recently in \cite{KH}. Except the Inoe
surfaces modeled on $Sol^4_0$, all other compact four-dimensional
solvmanifolds admitting a complex structure support also an
hyper-paracomplex structure  due to the results in
\cite{Pet,Kam,FPed} and Theorem~\ref{inoe}. It is also shown in
\cite{KH} that every complex structure on a compact 4-dimensional
solvmanifold is the canonical complex structure induced from the
left-invariant complex structure on the solvable Lie group $G$.
The four-dimensional Lie algebras admitting hyper-paracomplex
structure are classified in \cite{BVuk0}. A glance on Lie algebras
listed in \cite{BVuk0} leads to the conclusion that the Inoe
surfaces modeled on $Sol^4_0$ do not admit a hyper-paracomplex
structure induced from a left-invariant hyper-paracomplex
structure on the solvable Lie group $Sol^4_0$.

In view of Theorem~\ref{inq} and Theorem~\ref{inoe}, a naturally
arising question is whether the existence of a self-dual neutral
metric distinguishes the Inoe surfaces modeled on $Sol^4_1$ and
the Inoe surfaces modeled on $Sol^4_0$, i.e. whether there exists
a hyper-paracomplex structure on the Inoe surfaces modeled on
$Sol^4_0$.

We construct an anti-self-dual neutral metric which is not locally
conformally hyper-paraK\"ahler (hypersymplectic). We adapt the
Ashtekar at all \cite{Ash} formulation of the self-duality
Einstein equations to the case of neutral metric and modify the
Joyce's construction \cite{Joy} of hyper-complex structure from
holomorphic functions to get hyper-paracomplex structure.

Some properties of hyper-paracomplex and hyper-parahermitian
structures in higher dimensions are treated in \cite{IT,ITZ}.

We prove that the Nijenhuis tensor of a  Nearly paraK\"ahler
manifold is parallel with respect to the canonical connection. In
dimension six, we show that these spaces are Einsteinian but the
Ricci-flat case can not be excluded. This is in contrast with the
case of Nearly K\"ahler 6-manifolds which are Einsteinian with
positive scalar curvature. We involve twistor machinery to obtain
examples of Nearly paraK\"ahler manifolds. We adapt Salamon's
twistor construction on quaternionic manifold \cite{S1,S2,S3} to
the paraquaternionic situation. We consider the reflector space of
a paraquaternionic manifold as a higher dimensional analogue of
the reflector space of a 4-dimensional manifold with a metric of
neutral signature described in \cite{JR}. We show that the
reflector space of an Einstein self-dual non-Ricci flat 4 manifold
as well as the reflector space of a paraquaternionic K\"ahler
manifold admit both Nearly paraK\"ahler and almost paraK\"ahler
structures. We present homogeneous as well as non locally
homogeneous examples of 6-dimensional almost paraK\"ahler  and
Nearly paraK\"ahler  manifolds. However, all our examples of
Nearly paraK\"ahler  6-manifolds are Einstein spaces with non-zero
scalar curvature. To the best of the author's knowledge there are
no known examples of Ricci flat 6-dimensional Nearly paraK\"ahler
manifolds.

\begin{ack}
The final part of this paper was done during the visit of S.I. at
the Abdus Salam International Centre for Theoretical Physics,
Trieste, Italy. S.I. thanks the Abdus Salam ICTP for providing
support and an excellent research environment. S.I. is a  member
of the EDGE, Research Training Network HPRN-CT-2000-00101,
supported by the European Human Potential Programme. The research
is partially supported by Contract MM 809/1998 with the Ministry
of Science and Education of Bulgaria, Contract 586/2002 with the
University of Sofia "St. Kl. Ohridski". We are very grateful to V.
Apostolov for his comments and remarks. We also thank to R.
Cleyton for reading the manuscript and for comments.
\end{ack}

\section{Preliminaries}

Let $V$ be a real vector space of even dimension $2n$. An
endomorphism $P:V\rightarrow V$ is called a {\it paracomplex
structure} on $V$ if $P^2=1$ and the eigenspaces $V^{+1},V^{-1}$
corresponding to the eigenvalues $1$ and $-1$, respectively are of
the same dimension $n$, $V=V^+\oplus V^-$. Consider the algebra
$${\mathbb A} = \{x+\epsilon y, x,y \in {\mathbb R},
\epsilon^2=1\}$$ of paracomplex numbers over ${\mathbb R}$. As in
the ordinary complex case, ${\mathbb A^n}$ is identified with
$({\mathbb R^{2n}}, P)$, where $Pv= \epsilon v$. $P$ is called
{\it the canonical paracomplex structure} on ${\mathbb R^{2n}}$.

The notions of (almost) paracomplex, paraHermitian,
para-holomorphic, etc., objects are defined in the usual way over
the paracomplex numbers ${\mathbb A}$, instead of the complex
numbers ${\mathbb C}$. A survey on paracomplex geometry is
presented in \cite{CFG}.

A (1,1)-tensor filed $P$ on an 2n-dimensional smooth manifold $M$
is said to be an {\it almost product structure} if $P^2=1$. In
this case the pair $(M,P)$ is called {\it almost product
manifold}. An {\it almost paracomplex manifold} is an almost
product manifold $(M,P)$ such that the two eigenbundles $T^+M$ and
$T^-M$ associated with the two eigenvalues $\pm 1$ of $P$ have the
same rank. Equivalently, a splitting of the tangent bundle
$TM=TM^+\oplus TM^-$ of the subbundles $TM^{\pm}$ of the same
fiber dimension is called an almost paracomplex structure. A
smooth section of $TM^+$ is called \emph{(1,0)-vector field} while
a smooth section of $TM^-$ is said to be \emph{(0,1)-vector field}
with respect to the almost paracomplex structure. Such a structure
my alternatively be defined as a $G$-structure on $M$ with
structure group $GL(n,{\mathbb R})\times GL(n,{\mathbb R})$.

The Nijenhuis tensor $N$ of $P$ is defined by \cite{Yano}
\begin{equation*}
4N(X,Y) = [PX,PY]+[X,Y]-P[PX,Y]-P[X,PY].
\end{equation*}
The structure $P$ is said to be \emph{paracomplex} if $N=0$ \cite{Lib} which is
equivalent to the distributions on $M$ defined by $TM^{\pm}$ to be both completely
integrable \cite{KK}. The paracomplex manifold can also be characterized by
the existence of an atlas with paraholomorphic coordinate maps i.e. the
coordinate maps satisfying the para-Cauchy-Riemann equations \cite{Lib}
(see also \cite{KK}).

{\it An almost paraHermitian manifold} $(M,P,g)$ is a smooth
manifold endowed with an almost paracomplex structure $P$ and a
pseudo-Riemannian metric $g$ compatible in the sense that
\begin{equation*}
g(PX,Y)+g(X,PY)=0.
\end{equation*}
It follows that the metric $g$ is {\it neutral}, i.e. it has
signature (n,n) and the eigenbundles $TM^{\pm}$ are totally
isotropic with respect to $g$. Equivalently, an almost
paraHermitian manifold is a smooth manifold whose structure group
can be reduced to the real representation of the para-unitary
group

$$U(n,{\mathbb A})\cong \left\{\left(\begin{array}{c c}  A& 0\\
&\\ 0& (A^{-1})^t
\end{array}\right), A\in GL(n,\mathbb R)\right\}$$

isomorphic to $GL(n,\mathbb R)$.

Let $e_1,...,e_n,e_{n+1}=Pe_1,...,e_{2n}=Pe_n$ be an orthonormal
basis and denote $\epsilon_i=sign (g(e_i,e_i))=\pm 1,\quad
\epsilon_i=1, i=1,...,n,\quad \epsilon_j = -1, j=n+1,...,2n$.

The fundamental 2-form $F$ of an almost paraHermitian manifold is
defined by $$F(X,Y)=g(X,PY).$$ The covariant derivative of $F$
with respect to the Levi-Civita connection $\LC$ is expressed in
terms of $dF$ and $N$ in the following way (see e.g. \cite{KK})
\begin{gather}\label{df1}
2(\nabla^gF)(X;Y,Z)=-2g((\nabla^g_XP)Y,Z)=\\ \nonumber
dF(X,Y,Z)+dF(X,PY,PZ)+4N(PX;Y,Z).
\end{gather}
The Lee form $\theta$ is defined by $\theta = \delta F\circ P$,
where $\delta = -*d*$ is the co-differential with respect to $g$.
For 1-form $\alpha$ we use the notation $P\alpha(X)=-\alpha(PX)$.
Thus, $\theta=P\delta F$. We also have $$
\theta(X)=\sum_{i=1}^{2n}\epsilon_i(\nabla^gF)(e_i;e_i,PX)=
\frac{1}{2}\sum_{i=1}^{2n}dF(e_i,Pe_i,X)
=\sum_{i=1}^ndF(e_i,Pe_i,X).$$

Almost paraHermitian manifolds are classified with respect to the
decomposition in invariant and irreducible subspaces, under the
action of the structural group $U(n,{\mathbb A})$, of the vector
space of tensors satisfying the same symmetries as $\LC F$
\cite{Bej,GM}. We recall the defining conditions of some of the
classes:
\begin{enumerate}
\item[-] $\nabla^gF=0 \Leftrightarrow dF=0$, para-K\"ahler
manifolds;
\item[-] $N=0  \Leftrightarrow (\nabla^g_{PX}P)PY+(\nabla^g_XP)Y=0$,
paraHermitian manifolds \cite{Ol};
\item[-] $(\nabla^g_XP)X=0$, Nearly paraK\"ahler manifolds;
\item[-] $dF=0$, almost paraK\"ahler manifolds;
\item[-] $dF=\theta\wedge F, \quad d\theta=0$, paraHermitian manifolds locally
conformally equivalent to paraK\"ahler spaces \cite{GM,BO}.
\end{enumerate}

Examples of almost paraHermitian manifolds including the
non-compact hyperbolic Hopf and hyperbolic Calabi-Eckmann
manifolds \cite{Bej1} are collected in \cite{CFG}. Another source
of examples comes from the $k$-symmetric spaces, i.e. homogeneous
spaces defined by a Lie group automorphism of order $k$
\cite{Bal}. Almost paraHermitian manifolds are also called
\emph{almost bi-Lagrangian} \cite{JR,Kath}. They arise in relation
with the existence of Killing spinors of an indefinite neutral
metric \cite{Kath}.

\section{ParaHermitian connections}
A linear connection $\nabla$ on an almost paraHermitian manifold
$(M,g,P)$ is said to be {\it paraHermitian connection}, if it
preserves the paraHermitian structure, i.e. $\nabla g=\nabla P=0$.

In this section we define canonical paraHermitian connections in a (formally) similar way
as it was done in \cite{Ga1} for an almost Hermitian manifold.

We start with type decomposition of an element $B\in
\Lambda^2(TM)$. Denote $g(X,B(Y,Z)):= B(X;Y,Z)$. Let
$Bi(B):\Lambda^2(TM)\rightarrow\Lambda^3$ be the Bianchi projector
$$3Bi(B)(X;Y,Z)=B(X;Y,Z)+B(Y;Z,X)+B(Z;X,Y).$$ Further, we say that
$B$ is
\begin{enumerate}
\item[-] of type (1,1) if B(PX,PY)=-B(X,Y);
\item[-] of type (0,2) if B(PX,Y)=-PB(X,Y);
\item[-] of type (2,0) if B(PX,Y)=PB(X,Y).
\end{enumerate}
We will denote the corresponding  type-subspaces by $\Lambda^{1,1}$,
$\Lambda^{0,2}$, $\Lambda^{2,0}$, respectively, such
that $B=B^{1,1}\oplus B^{0,2}\oplus B^{2,0}$. The projections
are given by
\begin{gather*}
B^{1,1}(X,Y)=\frac{1}{2}\left(B(X,Y)-B(PX,PY)\right),\\
B^{0,2}(X,Y)=\frac{1}{4}\left(B(X,Y)+B(PX,PY)-PB(PX,Y)-PB(X,PY)\right),\\
B^{2,0}(X,Y)=\frac{1}{4}\left(B(X,Y)+B(PX,PY)+PB(PX,Y)+PB(X,PY)\right)
\end{gather*}
We define an involution $In:\Lambda^2(TM)\rightarrow \Lambda^2(TM)$ by
$In(B)(X;Y,Z,)=B(X;PY,PZ)$.

We may consider a 3-form $\psi$ as a totally skew-symmetric section of
$\Lambda^2(TM)$. It thus admits two different type decomposition:
\begin{enumerate}
\item[1.] decomposition as a 3-form: $\psi=\psi^+\oplus \psi^-$, where $\psi^+$
denotes the (1,2)+(2,1)-part and $\psi^-$-the (3,0)+(0,3)-part of $\psi$
given by
\begin{gather*}
\psi^+(X,Y,Z)=\frac{1}{4}\left(3\psi(X,Y,Z)-\psi(X,PY,PZ)-\psi(PX,Y,PZ)-
\psi(PX,PY,Z)\right),\\
\psi^-(X,Y,Z)=\frac{1}{4}\left(\psi(X,Y,Z)+\psi(X,PY,PZ)+\psi(PX,Y,PZ)+
\psi(PX,PY,Z)\right).
\end{gather*}
\item[2.] A type decomposition as an element of $\Lambda^2(TM)$.
\end{enumerate}
The two decompositions are related by
$
\psi^-=\psi^{0,2},\quad \psi^+=\psi^{2,0}+\psi^{1,1}.
$

Let $\nabla$ be any paraHermitian connection. Then we have
\begin{equation}\label{pot}
g(\nabla_XY,Z)-g(\nabla^g_XY,Z)=A(X;Y,Z),
\end{equation}
where $A\in\Lambda^2(TM)$ since $\nabla g=0$.

The torsion of $\nabla , T(X,Y)=\nabla_XY-\nabla_YX-\nabla_{[X,Y]}\in \Lambda^2(TM)$
and
\begin{equation}\label{tor}
T=-A+3Bi(A), \quad A=-T+\frac{3}{2}Bi(T), \quad Bi(A)=\frac{1}{2}Bi(T).
\end{equation}
We determine $\nabla$ in terms of its torsion.

Denote $d^aF(X,Y,Z):=-dF(PX,PY,PZ)$
we obtain easily  the following
\begin{pro}\label{treb}
On an almost paraHermitian manifold we have:
\begin{enumerate}
\item[a)] The Nijenhuis tensor is of type (0,2).
In particular it is trace-free, $tr(N)=0$.\\
The skew-symmetric part of $N$ is given by
\begin{equation}\label{nu2}
Bi(N)=\frac{1}{3}(d^aF)^-;
\end{equation}
\item[b)] The component $(\nabla^gF)^{1,1}=0$.
\item[c)] The component $(\nabla^gF)^{0,2}$ is determined by $N$:
\begin{gather}\label{nu3}
(\nabla^gF)^{0,2}(X;Y,Z)=dF^-(X,Y,Z)+2N(PX;Y,Z)=\\ \nonumber
N(PX;Y,Z)-N(PY;Z,X)-N(PZ;X,Y).
\end{gather}
\item[d)] The component $(\nabla^gF)^{2,0}$ is determined by $dF^+$:
\begin{equation}\label{df}
(\nabla^gF)^{2,0}(X;Y,Z)=\frac{1}{2}\left(dF^+(X,Y,Z)+dF^+(X,PY,PZ)\right)
\end{equation}
\end{enumerate}
\end{pro}
We describe the paraHermitian connections in the next
\begin{thm}\label{ter1}
Let $\nabla$ be a paraHermitian connection. Then
\begin{equation}\label{her}
T^{0,2}=-N,\qquad Bi(T^{2,0})-Bi(T^{1,1})=-\frac{1}{3}(d^aF)^+
\end{equation}
\indent For any 3-form $\psi^+$ of type (1,2)+(2,1) and any section $B_b$ of
$\Lambda^{1,1}(TM)$ satisfying $Bi(B_b)=0$ there exists a unique paraHermitian
connection whose torsion $T$ is given by the formula
\begin{equation}\label{tor1}
T=-N-\frac{1}{8}(d^aF)^+ -\frac{3}{8}In(d^aF)^++\frac{9}{8}\psi^++\frac{3}{8}
In(\psi^+)+B_b.
\end{equation}
\indent The corresponding paraHermitian connection is then equal to $\nabla^g+A$, where
$A$ is obtained from $T$ by \eqref{tor}.
\end{thm}
\begin{proof}
Since $\nabla P=0$ we get the first equality in \eqref{her} by
straightforward calculations. We calculate
$T^{2,0}-T^{1,1}=N+In(T),\quad 3Bi(In(T))=-d^aF$.
Apply \eqref{nu2} to derive
$3\left(Bi(T^{2,0})-Bi(T^{1,1})\right)=-d^aF+(d^aF)^-=-(d^aF)^+$ which completes
the proof of \eqref{her}.

Denote by $\psi^+$ the (1,2)+(2,1)-form $Bi(T^{2,0})+Bi(T^{1,1})$
and use \eqref{her} to get

\begin{equation}\label{t8}
Bi(T^{2,0})=\frac{1}{2}\left(\psi^+-\frac{1}{3}(d^aF)^+\right), \quad
Bi(T^{1,1})=\frac{1}{2}\left(\psi^++\frac{1}{3}(d^aF)^+\right).
\end{equation}
A linear connection $\nabla$ preserves the almost paracomplex
structure if and only if $A$ satisfies
$A(X;PY,Z)+A(X;Y,PZ)=(\nabla^gF)(X;Y,Z)$. By means of \eqref{tor}
the last equality is equivalent to
\begin{gather}\label{t9}
-T(X;PY,Z)-T(X;Y,PZ) + \frac{3}{2}\left(Bi(T)(X;PY,Z)+Bi(T)(X;Y,PZ)
\right)=\\ \nonumber =(\nabla^gF)(X;Y,Z).
\end{gather}
The first consequence of \eqref{t8} and \eqref{t9} is that the
(1,1)-part of $T$ which satisfies the Bianchi identity is free,
denote it by $T_b^{1,1}=B_b$. Take the (0,2) and (2,0) parts of
\eqref{t9}, apply \eqref{nu3}, \eqref{df} and use \eqref{her},
\eqref{t8} to get formula \eqref{tor1}.
\end{proof}
\begin{co}\label{co1}
Let $(M,g,P)$ be a 2n-dimensional almost paraHermitian manifold.
There exists paraHermitian connection on $M$ with totally
skew-symmetric torsion if an only if the Nijenhuis tensor is
totally skew-symmetric. In this case the connection is unique and
the torsion $T$ is given by
\begin{equation}\label{bis}
T=(d^aF)^+ -N
\end{equation}
\end{co}
\begin{proof}
Assume $T$ is a 3-form. Then $N$ is a 3-form due to \eqref{her} and $B_b=0$.
We claim $\psi^+=(d^aF)^+$. Indeed, $\psi^+=\frac{3}{4}T +\frac{1}{4}d^aF$.
On the other hand, $
\psi^+=Bi(T^{2,0})+Bi(T^{1,1}) = T+N=T+\frac{1}{3}(d^aF)^-$. Hence,
the claim follows. Substituting $\psi^+=(d^aF)^+$
into \eqref{tor1} we get \eqref{bis}. The corollary follows from
Theorem~\ref{ter1}
\end{proof}
We shall call this connection \emph{the Bismut connection}.
\begin{defn}
A paraHermitian connection is called
canonical if its torsion $T$ satisfies the following conditions
\begin{equation}\label{can}
T_b^{1,1}=0, \quad (Bi(T))^+= -\frac{2t-1}{3}(d^aF)^+
\end{equation}
for some real parameter $t$. We denote the corresponding connection
by $\nabla^t$.
\end{defn}
Combining \eqref{tor1} with \eqref{can} we get that the torsion
$T^t$ of $\nabla^t$ is given by
$$
T^t=-N-\frac{3t-1}{4}(d^aF)^+-\frac{t+1}{4}In(d^aF)^+.
$$
Any canonical connection is connected with the Levi-Civita connection
by
\begin{gather}\label{diff}
g(\nabla^t_XY,Z)=g(\nabla^g_XY,Z)-\frac{1}{2}g(\nabla^g_XP)(PY,Z)-\\ \nonumber
\frac{t}{4}\left((d^aF)^+(X,Y,Z)-(d^aF)^+(X,PY,PZ)\right).
\end{gather}
The paraHermitian connection with torsion 3-form is the canonical connection
given by $t=-1$. Another remarkable connection is the canonical connection
obtained for $t=0$ \cite{Yano1},
$$g(\nabla^0_XY,Z)=g(\nabla^g_XY,Z)-\frac{1}{2}g(\nabla^g_XP)(PY,Z), \quad
T^0=-N+\frac{1}{4}(d^aF)^+-\frac{1}{4}In(d^aF)^+.$$

Note that if $dF^+=0$ then the real line of the canonical connections
degenerates to a point $\nabla^0$ with torsion $T^0=-N$. Almost paraHermitian
manifolds satisfying the condition $dF^+=0$ are called \emph{quasi-paraK\"ahler}
or \emph{(1,2)-symplectic}. In view of Proposition~\ref{treb}, quasi-paraK\"ahler
manifolds are characterized by \cite{Yano1},
$(\nabla^g_{PX}F)(PY,Z)-((\nabla^g_XF)(Y,Z)=0$.

\subsection{Canonical connection on paraHermitian manifold}

We apply our previous discussion to a paraHermitian manifold, $N=0$.
\begin{thm}\label{cher}
Let $(M,g,P)$ be a 2n-dimensional paraHermitian manifold.
\begin{enumerate}
\item[a)] There exists a unique paraHermitian connection $\nabla^1$ on $M$
with torsion \\ $T^1\in \Lambda^{2,0}(TM)$ i.e. $T^1$ satisfies
$$T^1(PX,Y)=PT^1(X,Y).$$ This connection is the canonical
connection obtained by $t=1$ and given by
\begin{equation}\label{ch}
g(\nabla^1_XY,Z)=g(\nabla^g_XY,Z)-\frac{1}{2}dF(PX,Y,Z).
\end{equation}
\item[b)] The curvature $R^1:=[\nabla^1,\nabla^1]-\nabla^1_{[,]}$ is
of type (1,1) in the sense that
$$R^1(PX,PY)=-R^1(X,Y).$$
\end{enumerate}
\end{thm}
\begin{proof}
From $N=0$ we get $(d^aF)^-=0, d^aF=(d^aF)^+$.
Apply Theorem~\ref{ter1}. We have $B_b=0, \psi^+=Bi(T)=-\frac{1}{3}d^aF$
since $T^1\in \Lambda^{2,0}(TM)$. Hence, this is the canonical connection obtained for
 $t=1$ which proves a).

To prove b) we consider the paracomplex coordinate system
$(x^1,...,x^n,\bar{x}^1,...\bar{x}^n)$ around a point $p\in M$ such that
$\frac{\partial}{\partial x^1},...,\frac{\partial}{\partial x^n}$ is an
$+$-eigen-basis of $T_pM^+$ and $\frac{\partial}{\partial \bar{x}^1},...,
\frac{\partial}{\partial \bar{x}^n}$ is an
$-$-eigen-basis of $T_pM^-$, i.e. $P\frac{\partial}{\partial {x}^i}=
\frac{\partial}{\partial {x}^i}, \quad P\frac{\partial}{\partial \bar{x}^i}
=-\frac{\partial}{\partial \bar{x}^i}$. Then the metric and the fundamental
2-form are given by
$g=2g_{i\bar{j}}dx^idx^{\bar{j}}, \quad
F=F_{i\bar{j}}dx^i\wedge dx^{\bar{j}}, \quad F_{i\bar{j}}=-F_{\bar{j}i}=
-g_{i\bar{j}}.$

Summation in repeated indexes is always assumed. We use the
following convention: For a tensor $K$ of type (p,q),  the symbol
$\overline{K_{i_1,...\bar{i_q}}^{j_1,...,\bar{j_p}}}$ means
$K_{\bar{i_1},...i_q}^{\bar{j_1},...,j_p}$ .

We derive easily the expressions
$$dF_{ijk}=dF_{\bar i\bar j\bar k}=0,\quad dF_{ij\bar k}=
\frac{\partial g_{i\bar k}}{\partial x^j}-
\frac{\partial g_{j\bar k}}{\partial x^i},\quad
dF_{\bar i\bar jk}=\frac{\partial g_{k\bar j}}{\partial x^{\bar i}}-
\frac{\partial g_{k\bar i}}{\partial x^{\bar j}}=-\overline{dF_{ij\bar k}}.$$
Due to the  Koszul formula, the local components $\Gamma_{ij}^k$ of the Levi-Civita
connection are given by
\begin{gather}\label{gam}
\Gamma_{ij}^k=\frac{1}{2}g^{k\bar s}\left(
\frac{\partial g_{i\bar s}}{\partial x^j}+
\frac{\partial g_{j\bar s}}{\partial x^i}\right), \quad
\Gamma_{\bar i\bar j}^{\bar k}=\frac{1}{2}g^{\bar ks}\left(
\frac{\partial g_{\bar is}}{\partial x^{\bar j}}+
\frac{\partial g_{\bar js}}{\partial x^{\bar i}}\right),\\ \nonumber
\Gamma_{i\bar j}^k=\frac{1}{2}g^{k\bar s}dF_{\bar j\bar si}=
\Gamma_{\bar ji}^k, \quad
\Gamma_{i\bar j}^{\bar k}=\frac{1}{2}g^{s\bar k}dF_{si\bar j}=
\Gamma_{\bar ji}^{\bar k},\quad
\Gamma_{ij}^{\bar k}=\Gamma_{\bar i\bar j}^k=0.
\end{gather}
The local components
$C_{ij}^k$ of $\nabla^1$ are calculated from \eqref{ch} and \eqref{gam}
\begin{equation}\label{ch1}
C_{ij}^k=g^{k\bar s}\frac{\partial g_{j\bar s}}{\partial x^i},\quad
C_{\bar i\bar j}^{\bar k}=g^{s\bar k}\frac{\partial g_{s\bar j}}
{\partial x^{\bar i}}, \quad C_{i\bar j}^k=C_{\bar ij}^k=C_{ij}^{\bar k}=
C_{\bar i\bar j}^k=C_{i\bar j}^{\bar k}=C_{\bar ij}^{\bar k}=0
\end{equation}
The curvature tensor $R^1$ has the property $R^1\circ P=P\circ R^1$
since $\nabla^1P=0$. To prove b) it is sufficient to show $R^1_{ijk\bar l}=
R^1_{\bar i\bar j\bar kl}=0$ which is a direct consequence of
\eqref{ch1}.
\end{proof}
Further we shall call $\nabla^1$ \emph{the Chern connection}. This connection
coincides with the canonical compatible connection of the tangent bundle
viewed as a paraHermitian, paraholomorphic bundle of rank $n$ defined in
\cite{erd}.
\begin{co}\label{chern}
The curvature $R^g$ of the Levi-Civita connection of a
paraHermitian manifold satisfies the identities
$$R^g_{ijkl}=R^g_{\bar i\bar j\bar k\bar l}=0,$$ equivalently
\begin{gather*}
R^g(X,Y,Z,V)+R^g(PX,PY,PZ,PV)+R^g(X,Y,PZ,PV)+R^g(X,PY,Z,PV)+\\R^g(X,PY,PZ,V)+
R^g(PX,PY,Z,V)+ R^g(PX,Y,PZ,V)+R^g(PX,Y,Z,PV)=0.
\end{gather*}
The curvature $R^1$ and the torsion $T^1$
of the  Chern connection are given by
$$
R^1_{i\bar jk\bar l}=-g_{s\bar l}\frac{\partial
C^s_{ik}}{\partial x^{\bar j}}=
-g_{s\bar l}\frac{\partial}{\partial x^{\bar j}}\left(g^{s\bar m}
\frac{\partial g_{k\bar m}}{\partial x^i}\right),
\quad T_{\bar kij}=dF_{ij\bar k}.
$$
\end{co}
\subsection{Ricci forms of the  canonical connections}
For a linear connection $\nabla$ with curvature tensor $R$ on an almost
paraHermitian manifold of dimension $2n$ we have Ricci type tensors:
\begin{enumerate}
\item[-] the Ricci tensor $\rho(X,Y):=\sum_{i=1}^{2n}
\epsilon_iR(e_i,X,Y,e_i)$;
\item[-] the *-Ricci tensor
$\rho^*(X,Y):=\sum_{i=1}^{2n}
\epsilon_iR(e_i,X,PY,Pe_i)$;
\item[-] the Ricci form $r(X,Y):=-\frac{1}{2}\sum_{i=1}^{2n}
\epsilon_iR(X,Y,e_i,Pe_i) =-\sum_{i=1}^nR(X,Y,e_i,Pe_i)$.
\end{enumerate}
The  scalar curvatures are defined to be the corresponding trace:
\begin{enumerate}
\item[-] the scalar curvature $s=tr_g\rho=\sum_{i=1}^{2n}\epsilon_i\rho(e_i,e_i),$,
\item[-] the $*$-scalar curvature $s^*=tr_g\rho^*=2\sum_{i=1}^{2n}\epsilon_i\rho(e_i,e_i)$,
\item[-] the trace of the Ricci form $\tau=tr_gr=\sum_{i,j=1}^nr(e_i,Pe_i)$.
\end{enumerate}
For the Levi-Civita connection, we have the properties (see \cite{Ol})\\
$\rho^{g*}(X,Y)=r^g(X,PY), \quad \rho^{g*}(X,Y)+\rho^{g*}(PY,PX)=0$ and
consequently, $s^{g*}=\tau^g$.

To find relations between the Ricci forms of the canonical Hermitian connection
we consider the paraholomorphic canonical bundle $\Lambda_{\mathbb A}^n(TM)$.
Any linear connection preserving the structure $P$, i.e. preserving the
eigensubbubdles $TM^+$ and $TM^-$, induces a connection on the line bundle
 $\Lambda_{\mathbb A}^n(TM)$ with curvature equal to (-) its Ricci form. Let $s$ be a
 section of  $\Lambda_{\mathbb A}^n(TM)$. From \eqref{diff} we infer that
 $\nabla^ts=\nabla^0s+\frac{t}{2}P\theta\otimes s$. Consequently
$
 r^t=r^0-\frac{t}{2}d(P\theta)
$.
In particular the Ricci forms of the Bismut and Chern connection are related by
\begin{equation}\label{rform}
 r^{-1}=r^1+d(P\theta).
 \end{equation}
\section{ParaHermitian 4-manifold}

In this section we find a paraHermitian analogue of the Apostolov-Gauduchon
generalization \cite{AG1} of the Goldberg-Sachs theorem in General Relativity
(see e.g. \cite{PR}).
We prove the result using the properties of the Chern connection.

Let $(M,g)$ be an oriented pseudo-Riemannian 4-manifold with
neutral metric $g$ of signature $(+,+,-,-)$. This is equivalent,
on one hand to the existing of an almost paracomplex structure,
and on the other hand, to the existence of two kinds of  almost
complex structures. In a compact case the second property leads to
topological obstruction to the existence of neutral metric
expressed in terms of the signature and the Euler characteristic
\cite{Ma}.

The bundle $\Lambda^2M$ of real 2-forms of a neutral Riemannian 4-manifold
splits
\begin{equation}\label{4s1}
\Lambda^2M=\Lambda^+M\oplus \Lambda^-M,
\end{equation}
where $\Lambda^+M$, resp. $\Lambda^-M$ is the bundle of self-dual, resp. anti-self-dual
2-forms, i.e. the eigen-sub-bundle with respect to the eigenvalue $+1$, resp. $-1$,
of the Hodge $*$-operator acting as an involution on $\Lambda^2M$. We also
may consider the connected component $SO^+(2,2)$ of the structure group
$SO(2,2)$. This group has the splitting $SO^+(2,2)=SL(2)\times SL(2)$ which
defines two real vector bundles (of rank 2) $S^+$ and $S^-$ and $TM=S^+\oplus
S^-$ which induces the splitting \eqref{4s1}.

We will freely identify vectors and co-vectors via the metric $g$.

The self-dual part $W^+=\frac{1}{2}(W+*W)$ of the Weyl tensor $W$ is viewed as
a section of the bundle $\mathbb W^+ = Sym_0 \Lambda^+M$ of symmetric
traceless endomorphisms of $\Lambda^+M$.

Let $P$ be an almost paracomplex structure compatible with the metric $g$ such
that $(g,P)$ defines an almost paraHermitian structure. Then the fundamental
2-form $F$ is a section of $\Lambda^+M$ and has constant norm $2$. Conversely,
any smooth section of $\Lambda^+M$ with constant norm $2$ is the fundamental
2-form of an almost paracomplex structure. Our considerations in this section
are complementary to that in \cite{Apos1} in the sense that a section of
$\Lambda^+M$ with norm $-2$ can be considered as a K\"ahler form
of an almost complex structure, the case investigated in \cite{Apos1}.

We have
the following orthogonal splitting for $\Lambda^+M$
\begin{equation}\label{4s2}
\Lambda^+M=\mathbb R.F\oplus\Lambda_0^+M,
\end{equation}
where $\Lambda^+_0=\Lambda^{0,2}M\oplus \Lambda^{2,0}M$ denotes the bundle of
$P$-invariant real 2-forms $\phi$, $\phi(PX,PY)=\phi(X,Y)$.

In accordance with \eqref{4s2} the bundle $\mathbb W^+$ splits into three pieces as
follows:$$ \mathbb W^+=\mathbb W^+_1\oplus \mathbb W^+_2\oplus \mathbb W^+_3,$$
where
\begin{enumerate}

\item[-]
$\mathbb W^+_1=\mathbb M\times \mathbb R$ is the sub-bundle of elements
preserving \eqref{4s2} and acting by the homothety on the two factors,
hence the trivial line bundle generating by the elements
$\frac{3}{4}F\otimes F-\frac{1}{2}id;$
\item[-]
$\mathbb W^+_2$ is the sub-bundle of elements which exchange the two factors
in \eqref{4s2}: each element $\phi \in \Lambda^+_0M$ is identified with the element
$\frac{1}{2}(F\otimes\phi+\phi\otimes F)$;
\item[-]
$\mathbb W^+_3$ is the subbundle of elements preserving the splitting \eqref{4s2} and
acts trivially on the first factor $\mathbb R.F$, i.e. it is the space of those endomorphisms
of $\Lambda^+_0$ which are $P$-invariant.
\end{enumerate}
Thus, $W^+$ can be written in the form
\begin{equation}\label{4w1}
W^+=f\left(\frac{3}{4}F\otimes F-\frac{1}{2}id \right)+
\frac{1}{2}(F\otimes\phi+\phi\otimes F) + W^+_3,
\end{equation}
where $f$ is some real function.

In dimension $4$ the Lee form $\theta$ determines $dF$  completely by
\begin{equation}\label{4.1}
dF=\theta\wedge F.
\end{equation}
In particular, $d\theta$ is trace-free,
$\sum_{i=1}^2d\theta(e_i,Pe_i)=0$. Hence, the self-dual part $d\theta^+$ of
$d\theta$ is a section of $\Lambda^+_0M$.

To any 4-dimensional almost paraHermitian manifold with a Lee form $\theta$
one can associate the
canonical Weyl structure, i.e. a torsion-free connection $\nabla^w$
determined by the equation $\nabla^wg=\theta\otimes g$. The conformal scalar
curvature $k$ of an almost paraHermitian structure is defined to be the
scalar curvature of the canonical Weyl structure with respect to the metric $g$.
Then (see e.g. \cite{Gau})
\begin{equation}\label{ks1}
k=s-\frac{3}{2}\left(g(\theta,\theta)+2\delta\theta\right).
\end{equation}
The conformal scalar curvature is conformally invariant of weight $-2$, i.e. if $g'=f^{-2}g$
then $k'=f^2k$.
\subsection{Curvature of paraHermitian 4-manifold}
Let $(M,g,P)$ be a 4-dimensional paraHermitian manifold. The Chern connection
$\nabla^1$ and the Levi-Civita connection are related by $g(\nabla^1_XY,Z)=
g(\nabla^g_XY,Z) - \frac{1}{2}(\theta\wedge F)(PX,Y,Z)$ due to \eqref{ch}
and \eqref{4.1}. Consequently,
\begin{gather}\label{clc}
R^1(X,Y,Z,V)=R^g(X,Y,Z,V)-\frac{1}{2}d(P\theta)(X,Y)F(V,Z)+\\ \nonumber
\frac{1}{2}\left(L(Y,Z)g(V,X) -L(X,Z)g(V,Y)+L(X,V)g(Y,Z)-L(Y,V)g(Z,X) \right),
\end{gather}
where the tensor $L$ has the form
\begin{equation}\label{l1}
L(X,Y)=(\nabla^g_X\theta)Y+\frac{1}{2}\theta(X)\theta(Y)
-\frac{1}{4}g(\theta,\theta)g(X,Y).
\end{equation}
 The curvature $R^1$ is of type (1,1) according to Theorem~\ref{cher}. Then
\eqref{clc},\eqref{l1} imply, in local paraholomorphic coordinates, that
\begin{gather}\label{4r3}
R^g_{ijk\bar l}= -\frac{1}{2}\left(L_{jk}g_{i\bar l} -L_{ik}g_{j\bar l}
+d\theta_{ij}g_{k\bar l}\right),\quad R^g_{\bar i\bar j\bar kl}=
\overline{R^g_{ijk\bar l}};\\ \nonumber
R^g_{ij\bar k\bar l}= -\frac{1}{2}\left(L_{j\bar k}g_{i\bar l} -L_{i\bar k}
g_{j\bar l}-L_{j\bar l}g_{i\bar k} +L_{i\bar l}g_{j\bar k}\right), \quad
R^g_{\bar i\bar jkl}=\overline{R^g_{ij\bar k\bar l}},\\ \nonumber
 L_{ij}=\nabla^g_i\theta_j+\frac{1}{2}
\theta_i\theta_j, \qquad
L_{i\bar j}=\nabla^g_i\theta_{\bar j}+\frac{1}{2}\theta_i\theta_{\bar j}-
\frac{1}{4}|\theta|^2g_{i\bar j}.
\end{gather}
We take the traces in \eqref{clc}, \eqref{4r3}, \eqref{rform} and
use \eqref{ks1} to get  our technical
\begin{pro}\label{lemr}
The Ricci tensors and the scalar curvatures of a 4-dimensional paraHermitian manifold
satisfy the conditions
\begin{gather*}
\rho^g_{jk}=R^g_{ijk\bar i}+R^g_{ikj\bar i}=-\frac{1}{2}\left(
\nabla^g_j\theta_k + \nabla^g_k\theta_j+\theta_j\theta_k\right),\quad
\rho^g_{\bar j\bar k}=\overline{\rho^g_{jk}},
\\ \nonumber
\rho^{g*}_{jk}=-R^g_{ijk\bar i}+R^g_{ikj\bar i}=-\frac{1}{2}d\theta_{jk},
\quad \rho^{g*}_{\bar j\bar k}=\overline{\rho^{g*}_{jk}}
\\
\rho^g_{j\bar k}+\rho^{g*}_{j\bar k}=2R^g_{ij\bar k\bar i}=\left(
\frac{1}{2}\delta\theta
+\frac{1}{4}g(\theta,\theta)\right)g_{j\bar k},\\ \nonumber
s+s^*=2\delta\theta+g(\theta,\theta), \\ \nonumber k =-\frac{1}{2}(s+3s^*)=-\tau^{-1}.
\end{gather*}
In particular, the conformal scalar curvature is equal to $(-)$ the trace of the
Ricci form of the Bismut connection. Therefore, the trace of the Ricci form
of the Bismut connection is a conformal invariant of weight $-2$.
\end{pro}
We note that the expression of the (1,1)-part of the sum of the
two Ricci tensors and the formula for the sum of the two scalar
curvatures in Proposition~\ref{lemr} were obtained in \cite{Ol}.

The structure of $W^+$ on a 4-dimensional paraHermitian manifold is similar to
that of the Hermitian manifold presented in \cite{AG1}. We described it in the
following
\begin{lem}\label{w1}
On a 4-dimensional paraHermitian manifold the third component $W^+_3$ of $W^+$
 vanishes identically and the positive Weyl tensor is given by
\begin{equation}
W^+=\frac{k}{8}F\otimes F - \frac{k}{12}id -\frac{1}{4}\psi\otimes F-
\frac{1}{4}F\otimes\psi,
\end{equation}
where the two form $\psi$ is determined by the self-dual part of $d\theta_+,
\psi_{ij}=d\theta_{ij}$.
\end{lem}
\begin{proof}
On a 4-dimensional pseudo-Riemannian manifold the Weyl
tensor is expressed in terms of the normalized Ricci tensor
$h=-\frac{1}{2}\left(\rho -\frac{s}{6}g\right)$ as follows
\begin{gather}\label{w3}
W(X,Y,Z,V)=R^g(X,Y,Z,V)  -  h(X,Z)g(Y,V)+ \\ \nonumber
h(Y,Z)g(X,V)- h(Y,V)g(X,Z)+  h(X,V)g(Y,Z).
\end{gather}
The condition $W^+_3=0$ is a consequence of \eqref{w3} and
Corollary~\ref{chern}, due to the relation
$W_{ijkl}=R^g_{ijkl}=0$. According to \eqref{4w1} we have
$W(F)=W^+(F)=fF+\psi$. We calculate from \eqref{w3} applying
Proposition~\ref{lemr} that $$ W_{ijk\bar k}=-
\frac{1}{2}d\theta_{ij},\quad W_{i\bar jk\bar k}=
\frac{k}{6}g_{i\bar j}.$$ Hence, the lemma follows.
\end{proof}
Another glance at \eqref{rform} leads to the expression
$r^{-1}_{ij}=-d\theta_{ij}, \quad r^{-1}_{\bar i\bar
j}=d\theta_{\bar i\bar j}$ since the Ricci form of the Chern
connection is of type (1,1). The last equalities and
Lemma~\ref{w1} imply
\begin{pro}\label{bis1}
A 4-dimensional paraHermitian manifold is anti-self-dual ($W^+=0$) if and
only if the Ricci form of the Bismut connection is an anti-self-dual 2-form.
\end{pro}
Consider the co-differential of the positive Weyl tensor $\delta
W^+$ as an element of $\Lambda_0^2(T^*M)$. Then we have the
splitting $$\delta W^+=(\delta W^+)^+\oplus (\delta W^+)^-,$$
where $(\delta W^+)^+$ is a section of
$\Lambda_0^{(2,0)+(1,1)}(T^*M)$ while $(\delta W^+)^-$ is a
section of $\Lambda_0^{0,2}(T^*M)$. In particular, $(\delta
W^+)^-=0$ if and only if the co-differential of the whole Weyl
tensor vanishes on any three (1,0)-vectors.
\subsection{Proof of
Theorem~\ref{th4}} The equivalence a) $\Leftrightarrow$ b) is
proved in  Lemma~\ref{w1}.

The second Bianchi identity reads as $$\delta
W(X;Y,Z)=(\nabla^g_Yh)(Z,X)-(\nabla^g_Zh)(Y,X)$$ On (1,0) vectors
it  gives due to \eqref{gam} that
\begin{equation}\label{bia}
(\delta W^+)(X^{1,0};Y^{1,0},Z^{1,0})=
(\nabla^g_{Y^{1,0}}\rho)(Z^{1,0},X^{1,0})-(\nabla^g_{Z^{1,0}}\rho)(Y^{1,0},X^{1,0}).
\end{equation}
Assume $d\theta_{ij}=0$. Then  Proposition~\ref{lemr}, the Ricci
identities and  \eqref{4r3} imply
$\nabla^g_i\rho_{jk}-\nabla^g_j\rho_{ik}=0$. Hence, $(\delta
W^+)^- =0$ due to \eqref{bia}. The implication a) $\Rightarrow$ c)
is proved.

Let $(\delta W^+)^- =0$.  The local components of the Chern
connection and its torsion tensor are given by
$$C^k_{ij}=\Gamma^k_{ij}+\frac{1}{2}(\theta_i\delta_j^k-\theta_j\delta_i^k),
\quad T_{ij}^k=\theta_i\delta_j^k-\theta_j\delta_i^k.$$ The
equation \eqref{bia}, in terms of the Chern connection, takes the
form
\begin{equation}\label{biac}
\nabla^1_i\rho_{jk}-\nabla^1_j\rho_{ik}=\frac{3}{2}\left(\theta_j\rho_{ik}-
\theta_i\rho_{jk}\right).
\end{equation}
The Ricci identities for the Chern connection, $\nabla^1_i\nabla^1_j\theta_k-
\nabla^1_j\nabla^1_i\theta_k= \theta_j\nabla^1_i\theta_k -\theta_i\nabla^1_j\theta_k$,
 the first equality in Proposition~\ref{lemr} and \eqref{biac} yield
$$\nabla^1_id\theta_{jk}-\nabla^1_jd\theta_{ik}=\theta_kd\theta_{ij}-\frac{3}{2}
(\theta_id\theta_{jk}-\theta_jd\theta_{ik}).$$
Make a cyclic permutation in the latter then add the two and subtract the third of the
obtained equalities to get
\begin{equation}\label{cic1}
\nabla^1_id\theta_{jk}=-2\theta_id\theta_{jk}+
\frac{1}{2}(\theta_jd\theta_{ki}-\theta_kd\theta_{ji}).
\end{equation}
Take the covariant derivative in \eqref{cic1} and apply \eqref{cic1} to the obtained
result to derive
\begin{gather}\label{123}
\nabla^1_l\nabla^1_id\theta_{jk}=-2\nabla^1_l\theta_id\theta_{jk}+
\frac{1}{2}(\nabla^1_l\theta_jd\theta_{ki}-\nabla^1_l\theta_kd\theta_{ji})\\ \nonumber
+4\theta_i\theta_ld\theta_{jk}+\frac{5}{4}(\theta_i\theta_jd\theta_{lk} -
\theta_i\theta_kd\theta_{lj})-(\theta_j\theta_ld\theta_{ki}-\theta_k\theta_ld\theta_{ji})
\end{gather}
The Ricci identity $\nabla^1_i\nabla^1_jd\theta_{kl}-
\nabla^1_j\nabla^1_id\theta_{kl}=-\theta_i\nabla^1_jd\theta_{kl}+
\theta_j\nabla^1_id\theta_{kl}$ and \eqref{123} imply
\begin{gather*}
2d\theta_{li}d\theta_{jk}+\frac{3}{4}(\theta_i\theta_jd\theta_{kl}+
\theta_i\theta_kd\theta_{lj}-\theta_j\theta_ld\theta_{ki}-
\theta_k\theta_ld\theta_{ij})=\\
\frac{1}{2}(d\theta_{ij}\nabla^1_l\theta_k +d\theta_{ki}\nabla^1_l\theta_j-
d\theta_{kl}\nabla^1_i\theta_j-d\theta_{lj}\nabla^1_i\theta_k).
\end{gather*}
Change $l \leftrightarrow j, \quad i\leftrightarrow k$ into the
latter equality and sum up the results to obtain
$$4d\theta_{li}d\theta_{jk}=d\theta_{lk}d\theta_{ij}+d\theta_{lj}d\theta_{ki}.$$
From the last equality we easily infer
$5d\theta_{li}d\theta_{jk}=0$. Hence, $d\theta_{jk}=0$ which
completes the proof of Theorem~\ref{th4}.\hfill $\Box$\\ \indent
We consider the question of integrability of totally isotropic
real 2-plane supplementary
 distributions
on an oriented 4-dimensional neutral Riemannian manifold. Any such
splitting of the tangent bundle defines an almost paracomplex
structure compatible with the neutral metric, such that we get an
almost paraHermitian 4-manifold. The integrability of 2-plane
supplementary distributions is equivalent to the integrability of
the almost paracomplex structure. A necessary condition is the
vanishing of the third component $W^+_3$ of the positive Weyl
tensor, which is equivalent to the vanishing of the whole Weyl
tensor on the 2-plane distribution. Note that this is equivalent
to the vanishing of the whole curvature on the 2-plane i.e. the
identity in Corollary~\ref{chern} holds. This leads to fourth
order polynomial equation \cite{Akiv} (see also \cite{Apos1})
which can not have always real-root solutions. In the case of
existence, we give sufficient conditions for the integrability of
$P$ in the following
\begin{thm}\label{dis}
Let $(M,g)$ be an oriented neutral Riemannian 4-manifold with nowhere vanishing
positive Weyl tensor $W^+$. Suppose that $P$ is an almost paracomplex
structure such that  $W^+$ vanishes on each eigen-subbundle
determined by $P$, i.e. the component $W^+_3$ of $W^+$ with respect to $P$ vanishes.
Then any of the two following three conditions imply the third:
\begin{enumerate}
\item[i)] $W^+_2=0$;
\item[ii)] $(\delta W^+)^-=0$;
\item[iii)] the paracomplex  structure $P$ is integrable.
\end{enumerate}
\end{thm}
\begin{proof}
Observe that any smooth section $F$ of $\Lambda^+M$ with constant norm $2$ is
the fundamental 2-form of an almost paracomplex structure. Replacing $M$ by a two-fold
covering, if necessary, the positive Weyl tensor $W^+$ can be written in the form
\eqref{4w1}, where $f$ is a smooth function and $W^+_3=0$.

According to Theorem~\ref{th4} we have to show that i) and ii) imply iii).

Assume $W^+_2=0$. Then $W^+=\frac{3}{4}fF\otimes F-\frac{1}{2}f id$. Using the
definition of the Lee form, we calculate easily that
\begin{equation}\label{codif}
(\delta W^+)_X=\left(\frac{1}{2}Pdf(X)-\frac{3}{4}fP\theta(X)\right)F
-\frac{3}{4}f\nabla^g_{PX}F +\frac{1}{4}\left(df\wedge X+Pdf\wedge PX\right).
\end{equation}
The (0,2)-part of \eqref{codif} gives $0=(\delta W^+)^-=(\nabla^gF)^{0,2}$. Using
\eqref{nu3} we infer $N=0$.
\end{proof}
\begin{co}
 Let $(M,g,P)$ be an almost paraHermitian 4-manifold.
\begin{enumerate}
\item[i)] Suppose $W^+\not=0$ everywhere and $W^+_2=W^+_3=0$.
Then
$(\delta W^+)^+=0$ is equivalent to $d\left(|W^+|^{-\frac{2}{3}}F\right) =0.$
\item[ii)] Suppose $(M,g,P)$ is a paraHermitian 4-manifold. If it has
nowhere vanishing positive Weyl tensor then
$\delta W^+=0$ if and only if  $g'=|W^+|^{-\frac{2}{3}}g$ is a paraK\"ahler metric.\\
The Ricci tensor $\rho^g$ of $g$ is $P$-anti-invariant if and only if the
vector field $Pgrad_{g'}f$, where $f=|W^+|^{-\frac{1}{3}}$ is a Killing vector
field with respect to the paraK\"ahler metric $g'$.

In particular, a paraHermitian Einstein 4-manifold is either with everywhere
 vanishing positive Weyl tensor or  is globally conformal to a paraK\"ahler space.
 In the latter case there exists
non zero Killing vector field with respect to the paraK\"ahler metric.
\end{enumerate}
\end{co}
\begin{proof}
The (2,0)+(1,1)-part of \eqref{codif} yields
\begin{gather*}
(\delta W^+)^+_X=
\frac{3}{2}\left(\frac{1}{3}Pdf(X)-\frac{1}{2}fP\theta(X)
\right)F+\\
+\frac{3}{4}\left[\left(\frac{1}{3}df-\frac{1}{2}f\theta\right)\wedge
X+ \left(\frac{1}{3}Pdf-\frac{1}{2}fP\theta\right)\wedge
PX\right].
\end{gather*}
Assume $W^+_2=0$. Then the function $f$ is nowhere vanishing otherwise $W^+$ will
have zeros.
Moreover $|W^+|^2=(W^+(F,F))^2=4f^2$. The equation $(\delta W^+)^+=0$ is equivalent
to $\theta = \frac{2}{3}d ln f$. Thus we prove i).
The condition ii) is a consequence of i) and Theorem~\ref{dis}.
\end{proof}

\section{Nearly paraK\"ahler manifolds}

An almost paraHermitian manifold is called \emph{Nearly
paraK\"ahler (nearly bi-Lagrangian)} if the almost paraHermitian
structure is not para-K\"ahler and satisfies the identity
$$(\nabla^g_XP)X=0,\quad \Leftrightarrow (\nabla^g_XF)(Y,Z) +
(\nabla^g_YF)(X,Z)=0.$$ An example of nearly paraK\"ahlerian
6-manifold is given in \cite{Bej1}.

We denote the unique canonical connection $\nabla^0$ on a Nearly paraK\"ahler manifold
 by $\nabla$.

Applying the statements in Proposition~\ref{treb}, we get
\begin{pro}
A nearly paraK\"ahler manifold is quasi-K\"ahler, $dF^+=0$, the Nijenhuis tensor $N$
is a 3-form and the torsion $T$ of the unique canonical connection
is determined by the Nijenhuis tensor, $T=-N=P\nabla P$.
\end{pro}

Many properties of nearly paraK\"ahler manifolds are, in some
sense, formally very  similar to these of nearly K\"ahler
manifolds studied mainly by A.Gray \cite{Gr1,Gr2,Gr3}. Below we
follow roughly \cite{Gr3,Kir} (see also \cite{BM}).

\begin{pro}
On a nearly paraK\"ahler manifold the following identity holds
\begin{equation}\label{nk1}
R^g(X,Y,Z,V)+R^g(X,Y,PZ,PV)=g((\nabla^g_XP)Y,(\nabla^g_ZP)V).
\end{equation}
\end{pro}
\begin{proof}
The nearly paraK\"ahler condition implies $(\nabla^g_XP)(Y,PY)=0$.
Then, we get easily that
$R^g(X,Y,X,Y)+R^g(X,Y,PX,PY)=g((\nabla^g_XP)Y,(\nabla^g_XP)Y)$.
Polarizing the latter equality and using Bianchi identity, we
obtain \eqref{nk1}.
\end{proof}
Our crucial result in this section  is the following
\begin{thm}\label{par}
On a nearly paraK\"ahler manifold the Nijenhuis tensor is parallel with respect
to the canonical connection $\nabla$, $$\nabla N=-\nabla T=0.$$
\end{thm}
\begin{proof}
The curvature $R^g$ of the Levi-Civita connection and the curvature $R$ of the  canonical
 connection are related by
\begin{gather}\label{nkcur}
R^g(X,Y,Z,V)=R(X,Y,Z,V) - \frac{1}{2}(\nabla_XT)(Y,Z,V) +\frac{1}{2}
(\nabla_YT)(X,Z,V)\\ \nonumber
-\frac{1}{2}g\left(T(X,Y),T(Z,V)\right)
-\frac{1}{4}g\left(T(Y,Z),T(X,V)\right)
-\frac{1}{4}g\left(T(Z,X),T(Y,V)\right)
\end{gather}
Since $R\circ P=P\circ R$, the equation \eqref{nkcur} leads to
\begin{gather*}
R^g(X,Y,Z,V)+R^g(X,Y,PZ,PV)=\\
 -(\nabla_XT)(Y,Z,V)+(\nabla_YT)(X,Z,V)-g(T(X,Y),T(Z,V)).
 \end{gather*}
Comparing the latter equality with \eqref{nk1}, we derive
$(\nabla_XT)(Y,Z,V)-(\nabla_YT)(X,Z,V)=0$. Take the cyclic sum and add the result to conclude
 $\nabla T=0$.
\end{proof}
\begin{co}\label{nkcurv}
On a nearly paraK\"ahler manifold the following identities hold
\begin{gather}\nonumber
R^g(X,Y,Z,V)+R^g(X,Y,PZ,PV)+R^g(PX,Y,PZ,V)+R^g(PX,Y,Z,PV)=0, \\ \nonumber
R^g(X,Y,Z,V)=R^g(PX,PY,PZ,PV);\\ \nonumber
R(X,Y,Z,V)=R(Z,V,X,Y)=-R(PX,PY,Z,V)=-R(X,Y,PZ,PV);\\ \nonumber
\rho^g(PX,PY)=-\rho^g(X,Y), \quad \rho^{g*}(PX,PY)=-\rho^{g*}(X,Y)=-\rho^{g*}(Y,X),
\\ \nonumber
\rho(X,Y)=\rho(Y,X), \qquad r(PX,PY)=-r(X,Y),
\\ \label{nk0}
\rho^g(X,Y)-\rho(X,Y)=\frac{1}{2}\sum_{i=1}^ng(T(X,e_i),T(Y,e_i)), \\ \label{nk2}
\rho^g(X,Y)+\rho^{g*}(X,Y)=2\sum_{i=1}^ng(T(X,e_i),T(Y,e_i)),\\ \label{nk3}
\rho^{g*}(X,Y)=r^g(X,PY)=r(X,PY) -\frac{1}{2}\sum_{i=1}^ng(T(X,e_i),T(Y,e_i)),\\ \label{nk4}
3\rho^g(X,Y)-\rho^{g*}(X,Y)=4\rho(X,Y), \\ \label{nk5}
\rho^g(X,Y)+5\rho^{g*}(X,Y)=4r(X,PY).
\end{gather}
\end{co}
\begin{proof}
Put $\nabla T=0$ into  \eqref{nkcur} to get
\begin{gather}\label{nkcur1}
R^g(X,Y,Z,V)=R(X,Y,Z,V) \\ \nonumber
-\frac{1}{2}g\left(T(X,Y),T(Z,V)\right)
-\frac{1}{4}g\left(T(Y,Z),T(X,V)\right)
-\frac{1}{4}g\left(T(Z,X),T(Y,V)\right)
\end{gather}
All the identities in the corollary
are easy consequences of \eqref{nkcur1}
\end{proof}
\subsection{Nearly paraK\"ahler manifolds of dimension 6}
We recall that
a nearly paraK\"ahler manifold is said to be of
\emph{constant type} $\alpha\in \mathbb R$ if
\begin{equation}\label{contype}
g((\nabla^g_XP)Y,(\nabla^g_XP)Y)=
\alpha\left(g(X,X)g(Y,Y)-g^2(X,Y)+g^2(PX,Y)\right).
\end{equation}
In the Nearly K\"ahler case the constant type condition (with
positive constant $\alpha$) occurs only in dimension 6 and any
6-dimensional Nearly K\"ahler manifold is an Einstein manifold
with positive scalar curvature \cite{Gr3}. It is observed in
\cite{Kath} that in the Nearly paraK\"ahler case the constant type
phenomena occurs but the zero-value of $\alpha$ can not be
excluded. We describe the structure of the Ricci tensor in the
next
\begin{thm}\label{contyp}
Any 6-dimensional nearly paraK\"ahler manifold is an Einstein
manifold of constant type $\alpha\in \mathbb R$ and the following
relations hold
\begin{gather}
\rho^g=5\alpha g, \quad \rho^{g*}=-\alpha g,\quad \rho=4\alpha g.
\end{gather}
Consequently, the Riemannian scalar curvature $s^g=30\alpha$.

In particular, if $\alpha=0$ then the manifold is Ricci flat.
\end{thm}
\begin{proof}
Let $e_1,e_2,e_3,Pe_1,Pe_2,Pe_3$ be an orthonormal local basis of smooth vector fields.
The torsion $T$ of the canonical connection (or equivalently, the Nijenhuis tensor $N$)
is a 3-form of type (3,0)+(0,3). Therefore we may write $T(e_1,e_2)=ae_3+bPe_3$, where
$a$ and $b$ are smooth functions which turn to be constants because the torsion is
$\nabla$-parallel.
It is easy to calculate that \eqref{contype} holds with $\alpha = a^2-b^2$. Moreover,
we get the formula
\begin{equation}\label{kn1}
\sum_{i=1}^3g(T(X,e_i),T(Y,e_i))=2(a^2-b^2)g(X,Y)=2\alpha g(X,Y).
\end{equation}
The Nearly paraK\"ahler condition implies $(a,b)\not=(0,0)$. In
particular, the (3,0)+(0,3)-form $T$ is non-degenerate. On the
other hand, $\nabla T=0$, due to Theorem~\ref{par}. Hence, the
Ricci 2-form of the canonical connection vanishes as a curvature
of a flat line bundle. The condition $r=0$ and
Corollary~\ref{nkcurv} completes the proof.
\end{proof}
\begin{rem}
Using similar arguments as in the Nearly K\"ahler situation
\cite{Gr3} we derive that a Nearly paraK\"ahler manifold of
non-zero constant type has to be of dimension 6.
\end{rem}

\section{Examples, twistors and reflectors on paraquaternionic manifolds}

To obtain examples of Nearly paraK\"ahler manifolds we involve
twistor machinery. We are going to adapt Salamon's twistor
construction on quaternionic manifolds  \cite{S1,S2,S3} to the
paraquaternionic spaces.
\subsection{Paraquaternionic manifolds}
Both quaternions $H$ and paraquaternions $\tilde H$ are real Clifford algebras,
$H=C(2,0),\quad \tilde H=C(1,1)\cong C(0,2)$. In other words, the algebra
$\tilde H$ of paraquaternions is generated by the unity $1$ and the generators
$J_1,J_2,J_3$ satisfying the
\emph{paraquaternionic identities},
\begin{equation}\label{par1}
J_1^2=J_2^2=-J_3^2=1,\qquad J_1J_2=-J_2J_1=J_3.
\end{equation}

We recall the notion of almost paraquaternionic manifold
introduced by Libermann \cite{Lib}. An \emph{almost quaternionic
structure of the second kind}  on a smooth manifold consists of
two almost product structures $J_1,J_2$ and an almost complex
structure $J_3$, which mutually anti-commute, i.e. these
structures satisfy the  paraquaternionic identities \eqref{par1}.
Such a structure is also called \emph{complex product structure}
\cite{ASal,An}.

An \emph{almost hyper-paracomplex structure} on a 4n-dimensional
manifold $M$ is a triple $\tilde H=(J_{\alpha}), \alpha=1,2,3$,
where $J_a$,$a = 1,2$ are almost paracomplex structures
$J_a:TM\rightarrow TM$, and $J_{3}:TM\rightarrow TM$ is an almost
complex structure, satisfying the paraquaternionic identities
\eqref{par1}. When each $J_{\alpha},\alpha=1,2,3$ is an integrable
structure, $\tilde H$ is said to be a \emph{hyper-paracomplex
structure} on $M$. Such a structure is also called sometimes
\emph{pseudo-hyper-complex} \cite{D1}. Any hyper-paracomplex
structure admits a  unique torsion-free connection $\nabla^{ob}$
preserving $J_1,J_2,J_3$ \cite{ASal,An} called \emph{the complex
product connection}.

In fact an almost hyper-paracomplex structure is hyper-paracomplex
if and only if any two of the three structures $J_{\alpha}$ are
integrable, due to the following
\begin{pro}\label{nuij}
The Nijenhuis tensors $N_{\alpha}$ of an almost hyper-paracomplex
structure $\tilde H=(J_{\alpha}), \alpha=1,2,3$ are related by:
$$2N_{\alpha}(X,Y)=N_{\beta}(J_{\gamma}X,J_{\gamma}Y)-J_{\gamma}N_{\beta}(J_{\gamma}X,Y)-
J_{\gamma}N_{\beta}(X,J_{\gamma}Y)-J_{\gamma}^{2}N_{\beta}(X,Y)+$$
$$N_{\gamma}(J_{\beta}X,J_{\beta}Y)-J_{\beta}N_{\gamma}(J_{\beta}X,Y)-
J_{\beta}N_{\gamma}(X,J_{\beta}Y)-J_{\beta}^{2}N_{\gamma}(X,Y)$$
\end{pro}
\begin{proof}
The formula follows by very definitions with long but standard computations.
\end{proof}
We note that during the preparation of the manuscript the formula in
the Proposition~\ref{nuij} appeared in the context of Lie algebras in \cite{BVuk0}.

An \emph{almost paraquaternionic structure} on $M$ is a rank-3
subbundle $P \subset End(TM)$ which is locally spanned  by an
almost hyper-paracomplex structure $\tilde H=(J_{\alpha})$; such a
locally defined triple $\tilde H$ will be called  admissible basis
of $P$. A linear connection $D$ on $TM$ is called
\emph{paraquaternionic connection} if $D$ preserves $P$, i.e.
there exist locally defined 1-forms $\omega_{\alpha}, \alpha
=1,2,3$ such that
\begin{equation}\label{zzv}
DJ_1=-\omega_3\otimes J_2+\omega_2\otimes J_3, \quad
DJ_2=\omega_3\otimes J_1+\omega_1\otimes J_3, \quad
DJ_3=\omega_2\otimes J_1+\omega_1\otimes J_2.
\end{equation}
Consequently, the curvature $R^D$ of $D$ satisfies the relations
\begin{gather}\label{rel1}\nonumber
[R^D,J_1]=-A_3\otimes J_2+A_2\otimes J_3, \\ [R^D,J_2]=A_3\otimes
J_1+A_1\otimes J_3, \\ \nonumber [R^D,J_3]=A_2\otimes J_1+A\otimes
J_2,\\ \nonumber A_1=d\omega_1+\omega_2\wedge\omega_3, \quad
A_2=d\omega_2+\omega_3\wedge\omega_1, \quad
A_3=d\omega_3-\omega_1\wedge\omega_2.
\end{gather}

An almost paraquaternionic structure is said to be a
\emph{paraquaternionic} if there is a torsion-free paraquaternionic connection.

A \emph{$P$-Hermitian metric} is a pseudo Riemannian metric which is compatible
with the (almost) hyper-paracomplex structure $\tilde H=(J_{\alpha}), \alpha=1,2,3$
in the sense
that the metric $g$ is skew-symmetric with respect  to each
$J_{\alpha}, \alpha=1,2,3$, i.e.
\begin{equation}\label{par2}
g(J_1.,J_1.)=g(J_2.,J_2.)=-g(J_3.,J_3.)=-g(.,.).
\end{equation}
The metric $g$ is necessarily of neutral signature (2n,2n). Such a structure is
called \emph{(almost) hyper-paraHermitian structure}.

An almost paraquaternionic (resp. paraquaternionic) manifold with
P-Hermitian metric is called an \emph{almost paraquaternionic
Hermitian} (resp. \emph{paraquaternionic Hermitian}) manifold. If
the Levi-Civita connection of a paraquaternionic Hermitian
manifold is paraquaternionic connection, then the manifold is said
to be \emph{paraquaternionic K\"ahler} manifold. This condition is
equivalent to the statement that the holonomy group of $g$ is
contained in $Sp(n,\mathbb R)Sp(1,\mathbb R)$ for $n\ge 2$
\cite{GRio,Vuk}. A typical example is the paraquaternionic
projective space endowed with the standard paraquaternionic
K\"ahler structure \cite{BL1}. Any paraquaternionic K\"ahler
manifold of dimension $4n\ge 8$ is known to be Einstein with
scalar curvature $s$ \cite{GRio,Vuk}.  If on a paraquaternionic
K\"ahler manifold there exists an admissible basis $(\tilde H)$
such that each $J_{\alpha}, \alpha=1,2,3$ is parallel with respect
to the Levi-Civita connection, then the manifold is said to be
\emph{hyper-paraK\"ahler}. Such manifolds are also called
\emph{hypersymplectic} \cite{Hit}, \emph{neutral hyper-K\"ahler}
\cite{Kam,FPed}. The equivalent characterization is that the
holonomy group of $g$ is contained in $Sp(n,\mathbb R)$ if $ n\ge
2$ \cite{Vuk}.

When $n\ge 2$, the paraquaternionic condition, i.e. the existence of
torsion-free paraquaternionic connection is a strong condition
which is equivalent to the 1-integrability of the associated
$GL(n,\tilde H)Sp(1,\mathbb R)\cong GL(2n,\mathbb R)Sp(1,\mathbb R)$-
structure \cite{An,ASal}. Such a structure is a type of a para-conformal structure
\cite{BE} as well as a type of generalized hypercomplex structure \cite{Biel}.

\subsection{Hyper-paracomplex structures on 4-manifold}
For $n=1$ an almost paraquaternionic structure is the same as an
oriented neutral conformal structure and turns out to be always
paraquaternionic \cite{D1,GRio,Vuk,BVuk0}.
 The existence of a (local) hyper-paracomplex
structure is a strong condition because of the next
\begin{thm}\label{inq}
If on a 4-manifold there exists a (local) hyper-paracomplex structure
then the corresponding neutral conformal structure is anti-self-dual.
\end{thm}
\begin{proof}
Let $(g,(J_{\alpha}),\alpha=1,2,3)$ be an almost
hyper-paraHermitian structure with fundamental 2-form $F_{\alpha}$
associated to each $J_{\alpha}$.  Denote by
$\theta_1,\theta_2,\theta_3$ the corresponding Lee forms (defined
by $\theta_{\alpha}=-\delta F_{\alpha}\circ J_{\alpha}^3$).
\begin{lem}\label{mlq}
The
 structure $(g,(J_{\alpha}),\alpha=1,2,3)$
is a hyper-paracomplex structure, if and only if the
three Lee forms coincide, $\theta_1=\theta_2=\theta_3$.
\end{lem}
\begin{proof}
The Levi-Civita connection satisfies \eqref{zzv}. Consequently the Nijenhuis tensors
obey
\begin{equation}\label{nnn}
N_{\alpha}=-B_{\alpha}\otimes J_{\beta}+ J_{\beta}\otimes  B_{\alpha}-
J_{\alpha}B_{\alpha}\otimes J_{\gamma}+ J_{\gamma}\otimes J_{\alpha}B_{\alpha},\quad
B_{\alpha}=\omega_{\beta}-J_{\alpha}^3\omega_{\gamma}.
\end{equation}
Simple calculations using \eqref{zzv} give

\centerline{$\theta_1=-J_2\omega_2+J_3\omega_3, \quad
\theta_2=J_1\omega_1+J_3\omega_3, \quad
\theta_3=-J_2\omega_2+J_1\omega_1$.}
The last three identities and \eqref{nnn} yield

\centerline{$J_1(\theta_2-\theta_1)=B_3,\quad
J_2(\theta_2-\theta_3)=B_1,\quad
J_3(\theta_3-\theta_1)=B_2$.}
Another glance at \eqref{nnn} completes the proof of the lemma.
\end{proof}
Suppose that each $J_{\alpha},\alpha=1,2,3$ is integrable. Denote
the common Lee form by $\theta$ and take the 3-form $T$ to be the
Hodge-dual to $\theta$ with respect to $g$. We have the identities
$T=*\theta=-\theta\circ J_1\wedge F_1= -\theta\circ J_2\wedge F_2=
+\theta\circ J_3\wedge F_3$. Then the Bismut connections of the
three structures coincide, i.e. the linear connection
$\nabla^b:=\nabla^g +\frac{1}{2}T$ preserves the metric and each
$J_{\alpha},\alpha =1,2,3$. Therefore, each fundamental two form
is parallel with respect to this connection,
$\nabla^bF_{\alpha}=0, \alpha=1,2,3$. Consequently, the 2-form
$\Phi_1 = F_2+F_3$ is $\nabla^b$-parallel, $\nabla^b\Phi_1=0$. The
2-form $\Phi_1$ is a (2,0)+(0,2)-form with respect to $J_1$.
Hence, the Ricci form of the Bismut connection vanishes and
$W^+=0$
 due to Proposition~\ref{bis1}
\end{proof}
Note that the integrability condition, Lemma~\ref{mlq}, in the
case of hyper-complex structure is due to F.Battaglia and
S.Salamon (see \cite{GT}).

The universal cover $\widetilde{SL(2,\mathbb R)}$ of the Lie group
${SL(2,\mathbb R)}$ admits a discrete subgroup  $\Gamma$ such that
the quotient space $\widetilde{(SL(2,\mathbb R)}/\Gamma)$ is a
compact 3-manifold  \cite{Mil,RV,Scot}. Such a space has to be
Seifert fibre space \cite{Scot} and all the quotients are
classified in \cite{RV}. The compact 4-manifold $M=S^1\times
\widetilde{(SL(2,\mathbb R)}/\Gamma)$ admits a complex structure
and is known as Kodaira-Thurston surface modeled on $S^1\times
\widetilde{SL(2,\mathbb R)}$ \cite{Wall}.

\subsection{Proof of Theorem~\ref{kod-t}} Let
$S^2_2(1)=\{\mathbb R^4\ni (a,b,c,d) : a^2+b^2-c^2-d^2=1\}$ be the
unit pseudo-sphere with respect to the standart neutral metric in
$\mathbb R^4$. We consider the so-called \emph{hyperbolic Hopf
manifold} $\mathbb R\times S^2_2(1)$ isomorphic to the Lie group
$\mathbb R\times SL(2,\mathbb R)$. The Lie algebra $\mathbb
R\times sl(2,\mathbb R)\cong gl(2,\mathbb R)$ has a basis
$\{W,X,Y,Z\}$ with $Z$ central and non-zero brackets given by \\
\centerline{$[X,Y]=W, \quad [Y,W]=-X, \quad [W,X]=Y$.} An almost
paracomplex structure on $\mathbb R\times SL(2,\mathbb R)$ is
constructed in \cite{Bej1}. The Lie algebra $\mathbb R\times
sl(2,\mathbb R)$ supports a hyper-paracomplex structure given by
\cite{ASal,BVuk0}\\ \centerline{$J_3Z=X, \quad J_3Y=W,\quad
J_2Z=Y,\quad J_2X=-W$.} We pick a compatible neutral metric $g$,
in the corresponding conformal class, defined such that  the basis
$\{W,X,Y,Z\}$ is an orthonormal basis, $X,Z$ have norm $1$ while
$Y,W$ have norm $-1$, $g(X,X)=g(Z,Z)=-g(W,W)=-g(Y,Y)=1$.
\begin{lem}\label{hopf}
The invariant hyper-paraHermitian structure on $\mathbb R\times
SL(2,\mathbb R)$, described above, is non-flat conformally
equivalent to a flat hyper-paraK\"ahler structure.

More precisely, the Lee form $\theta=-Z$ is $\nabla^g$-parallel
and the complex product  connection coincides with the Levi-Civita
connection of the flat hyper-paraK\"ahler
 metric $g^{ob}=e^{-t}g$, where $t$ is the
local coordinate on $\mathbb R$.
\end{lem}
\begin{proof}
The Koszul formula
gives the following non-zero terms:
\begin{gather*}
2\nabla^g_XY=W,\quad 2\nabla^g_YW=-X,\quad  2\nabla^g_WX=Y,\\
2\nabla^g_XW=-Y,\quad 2\nabla^g_YX=-W,\quad  2\nabla^g_WY=X.
\end{gather*}
It is easy to check that $g$ is not flat and  $\theta=-Z =-dt$
satisfies $\nabla^g\theta=0$. The Levi-Civita connection of the
conformal metric $g'=e^{-t}g$ is determined by\\
\centerline{$2\nabla^{g'}_AB:=2\nabla^g_AB-\theta(A)B-\theta(B)A+g(A,B)\theta$.}
It is straightforward to verify that $\nabla^{g'}$  preserves
$J_1,J_2,J_3$. Hence, it is the complex product connection and the
metric $g'$ is hyper-paraK\"ahler. It is not difficult to
calculate that the connection $\nabla^{g'}$ is flat which proves
the lemma.
\end{proof}
An (left) invariant Weyl-flat hyper-paraHermitian structure on
$\mathbb R\times SL(2,\mathbb R)$ is just the neutral product of
the standard Lorentz metric of constant sectional curvature on the
unit pseudo-sphere $S^2_2(1)$, induced by the neutral metric on
$\mathbb R^4$, and the flat metric on $\mathbb R$. In coordinates
$(x,y,z,t)$, it has the form

\centerline{$ds^2=(\cosh y)^2(\cosh z)^2dxdx+dtdt-(\cosh z)^2dydy-dzdz.$}

The lelft-invariant Weyl-flat hyper-paraHermitian structure on
$\mathbb R\times \widetilde{SL(2,\mathbb R)}$ described in
Lemma~\ref{hopf} descends to $M=S^1\times \widetilde{(SL(2,\mathbb
R)}/\Gamma)$. The descended structure is not globally conformal to
a hyper-paraK\"ahler structure since the closed Lee form  $\theta$ is
actually a 1 form on the circle $S^1$ and therefore can not be
exact. Hence, the proof of Theorem~\ref{kod-t} is completed.\hfill
$\Box$\\ \indent The 4-dimensional Lie algebras admitting a
hyper-paracomplex structure were classified recently in
\cite{BVuk0}. It is shown in \cite{BVuk0} that exactly 10 types of
Lie algebras admit a hyper-paracomplex structure.
Theorem~\ref{inq} tells us that the corresponding neutral metrics
are anti-self-dual. We show below that some of them are not
conformally flat.

Note that all 4-dimensional Lie groups admitting anti-self-dual
non Weyl-flat Riemannian metric are classified in \cite{DS}.

\begin{example}\label{nex1}
\begin{enumerate}
\item[] We recall the construction of hyper-paracomplex structures on some
4-dimensional Lie algebras keeping the notations in \cite{BVuk0}.

\item[i)] Consider the solvable Lie algebra PHC5 with a basis $\{X,Y,Z,W\}$, non-zero
bracket $[X,Y]=X$ and hyper-paracomplex structure given by

\centerline{$J_3Z=W,\quad J_3X=Y,\quad J_2Z=W,\quad J_2X=Y-Z,\quad J_2Y=X+W$.}

\noindent Consider the oriented basis $A=X,\quad B=Y,\quad C=Y-Z,\quad D=-X-W$ and
pick a compatible neutral metric $g$ with non-zero values on the basis
$\{A,B,C,D\}$ given by $g(A,A)=g(B,B)=-g(C,C)=-g(D,D)=1$. The metric $g$ on the
corresponding simply connected solvable Lie group is
conformally hyper-paraK\"ahler (hypersymplectic) since the Lee form $\theta=B-C$ is closed
and therefore exact.
 It is anti-self-dual
metric with non-zero Weyl tensor because its curvature
$R^g(A,B,C,D)=1$. In local coordinates $\{x,y,z,t\}$, the metric
is given by
$$ ds^2=e^{2y}dxdx+dydy-e^{-y}(dxdt +dtdx) + (dydz +
dzdy). $$
\item[ii)] Consider the solvable Lie algebras PHC6, PHC9, PHC10 defined by non-zero brackets:
\item[PHC6] $[X,Y]=Z, [X,W]=X+aY+bZ, [W,Y]=Y$
\item[PHC9] $[Z,W]=Z, [X,W]=cX+aY+bZ, [Y,W]=Y, c\not=0$
\item[PHC10] $[Y,X]=Z, [W,Z]=cZ, [W,X]=\frac{1}{2}X+aY+bZ, [W,Y]=(c-\frac{1}{2})Y, c\not=0$

These algebras admit a hyper-paracomplex structure defined by

\centerline{$J_3Z=Y,\quad J_3X=W,\quad J_2Z=Y,\quad J_2X=W-Z,\quad J_2W=X+Y$.}
\noindent Consider the oriented frame $A=X, B=W, C=W-Z, D=-X-Y$. A compatible metric $g$ is defined
such that the frame $\{A,B,C,D\}$ is orthonormal with $g(A,A)=g(B,B)=-g(C,C)=-g(D,D)=1$.
The Lee forms of these hyper-paraHermitian structures are closed and the curvature satisfies

\item[PHC6] $R^g(A,B,C,D)=(1-a);$
\item[PHC9]  $R^g(A,B,C,D)=\frac{1}{2}(2c^2-3c-2ac+2a+1);$
\item[PHC10] $R^g(A,B,C,D)=\frac{1}{2}(c^2+2ac-c);$

\noindent Clearly there are constants (a,b,c) such that the
corresponding Lie algebras admit anti-self-dual neutral metric
with non-zero Weyl tensor. For example, let us take $c=-2,\quad
a=b=0$ in the Lie algebra $PHC9$ described in Example~\ref{nex1},
ii). The Lee form $\theta=B-C$ is not $\nabla^g$-parallel but
closed and the Weyl curvature does not vanish because
$R(A,B,C,D)=15/2$. In coordinates $x,y,z,t$ the left invariant
vector fields $A,B,C,D$ can be expressed as follows

$A=e^{-2t}\frac{\partial}{\partial x}, \quad
B=\frac{\partial}{\partial t}, \quad C=\frac{\partial}{\partial
t}- e^t\frac{\partial}{\partial z},\quad
D=-e^{-2t}\frac{\partial}{\partial x} -e^t\frac{\partial}{\partial
y}.
$

The invariant neutral anti-self-dual metric with non-zero Weyl
tensor has the form

\centerline{$ds^2=e^{4t}dxdx+dtdt-e^t(dxdy+dydx)+e^{-t}(dtdz+dzdt)$.}
\end{enumerate}
\end{example}
It turns out that the conformal structure $[g]$ induced by the
invariant hyper-paracomplex structure on the corresponding simply
connected 4-dimensional Lie group is actually generated by a
hyper-paraK\"ahler (hypersymplectic) structure, since  the Lee
form $\theta$ is closed (and therefore exact) in all 10 possible
cases described in \cite{BVuk0}. On some of them the Lee form is
zero and the structure is hyper-paraK\"ahler (hypersymplectic).
\begin{example}
The solvable Lie group corresponding to the Lie algebra defined in
\cite{ASal}  and obtained from $PHC9$ for $c=-1, \quad a=b=0$,
 posses an invariant hyper-paraK\"ahler (hypersymplectic) structure with
non-zero Weyl tensor since the Lee form vanishes and the curvature
has non-zero value on an orthonormal basis.
\end{example}
 Summarizing, we get
\begin{pro}\label{qup}
Any one of the nine simply connected solvable Lie groups corresponding to a
solvable 4-dimensional
Lie algebra admitting hyper-paracomplex structure  supports a
hyper-paraK\"ahler (hypersymplectic) structure.
\end{pro}
\begin{rem}
The hyper-paraK\"ahler (hypersymplectic) structures on the nine
solvable Lie groups mentioned in Proposition~\ref{qup}, are not
left-invariant in general. There are left-invariant
hypersymplectic structure on exactly four cases
  according to the recent classification
of the  hypersymplectic 4-dimensional Lie algebras
 \cite{Andr}.
\end{rem}
Due to the Malcev theorem \cite{Mal}, the 4-dimensional nilpotent
Lie group $H$  has a discrete subgroup $\Gamma$ such that the
quotient $M=H/\Gamma$ is a compact nil-manifold, the Kodaira
surface. It is known that these surfaces admit a
hyper-paraK\"ahler (hypersymplectic) structure \cite{Kam}, see
also \cite{FPed}.

Consider the solvable Lie algebra $sol^4_1$ defined by non-zero
brackets: $$[X,Y]=Z, [X,W]=X, [W,Y]=Y.$$ This Lie algebra can be
obtained by taking $a=b=0$ in the Lie algebra $PHC6$ described in
Example~\ref{nex1}, ii).

The corresponding solvable Lie group is known to be $Sol_1^4$. The
geometric structures modeled on this group appear as one of the
possible geometric structures on 4-manifold \cite{Wall}. The
compact quotients of $Sol^4_1$ by a discrete group $\Gamma$
constitute the Inoe surfaces modeled on $Sol^4_1$ \cite{Wall}.
\subsection{Proof of Theorem~\ref{inoe}}
A  hyper-paracomplex structure on the Lie algebra $sol^4_1$ is
given by

\centerline{$J_3Z=Y,\quad J_3X=W,\quad J_2Z=Y,\quad J_2X=W-Z,\quad
J_2W=X+Y$.}

\noindent Consider the oriented frame $A=X, B=W, C=W-Z, D=-X-Y$. A
compatible metric $g$ is defined such that the frame $\{A,B,C,D\}$
is orthonormal with $g(A,A)=g(B,B)=-g(C,C)=-g(D,D)=1$.
 The Lee form $\theta=B-C$ is not
$\nabla^g$-parallel but closed and the Weyl curvature does not
vanish because $R(A,B,C,D)=1$.

In coordinates $x,y,z,t$, the left invariant vector fields
$A,B,C,D$ on $Sol^4_1$ can be expressed as follows

$A=e^{-t}\frac{\partial}{\partial x}, \quad
B=\frac{\partial}{\partial t}, \quad
C=\frac{\partial}{\partial t}-
\frac{\partial}{\partial z},\quad
D=-e^{-t}\frac{\partial}{\partial x}
-e^t\frac{\partial}{\partial y}-e^tx\frac{\partial}{\partial z}.
$

The left invariant neutral anti-self-dual metric with non-zero
Weyl tensor on $Sol^4_1$ has the form

\centerline{$ds^2=e^{2t}dxdx+dtdt-(dxdy+dydx)+
(dtdz+dzdt)-x(dtdy+dydt)$.}

The left invariant hyper-paracomplex structures on $Sol_1^4$,
described above, descends to the Inoe surfaces modeled on
$Sol^4_1$. Theorem~\ref{inq} completes the proof of
Theorem~\ref{inoe}. \hfill $\Box$\\ \indent All local
hyper-paracomplex structures on 4-manifold, we have presented, are
locally conformal to hypersymplectic structures. We shall
construct a local hyper-paracomplex structure which is not
conformally equivalent to a hyper-paraK\"ahler (hypersymplectic),
i.e. its Lee form $d\theta\not=0$. We adapt the Ashtekar at all
\cite{Ash} formulation of the self-duality Einstein equations to
the case of neutral metric and modify  the Joyce's construction
\cite{Joy} of hyper-complex structure from holomorphic functions.
\begin{example}
Let $V_1,V_2,V_3,V_4$ be a vector fields on an oriented 4-manifold
$M$ forming an oriented basis for $TM$ at each point. Then
$V_1,..., V_4$ define a neutral conformal structure $[g]$ on $M$.
Define an almost hyper-paracomplex structure $(J_2,J_3)$ by the
equations

\centerline{$J_3V_1=-V_2,\quad J_3V_3=V_4,\qquad J_2V_1=-V_4, \quad J_2V_2=V_3$.}
Suppose that $V_1,..., V_4$ satisfy the three vector field equations
\begin{equation}\label{joi1}
[V_1,V_2]+[V_3,V_4]=0,\quad [V_1,V_3]+[V_2,V_4]=0, \quad [V_1,V_4]-[V_2,V_3]=0.
\end{equation}
It is easy to check that these equations imply the integrability
of $(J_2,J_3)$, i.e. $(J_2,J_3)$ is a hyper-paracomplex structure
which is compatible with the neutral conformal structure $[g]$.
Hence, $[g]$ is anti-self-dual, due to Theorem~\ref{inq}.

The neutral Ashtekar at all equation \eqref{joi1} may be written
in a complex form
\begin{equation}\label{joy2}
[V_1+iV_2,V_1-iV_2]+[V_3+iV_4,V_3-iV_4]=0, \quad [V_1+iV_2,V_3-iV_4]=0.
\end{equation}
Let $M$ be a complex surface, let $(z^1,z^2)$ be local holomorphic coordinates, and
define  $V_1,..., V_4$ by
$$V_1+iV_2=f_1\frac{\partial}{\partial z^1} + f_2\frac{\partial}{\partial z^2},
\quad
V_3+iV_4=f_3\frac{\partial}{\partial z^1}+f_4\frac{\partial}{\partial z^2},$$
where $f_j$ is a complex function on $M$. Substituting into \eqref{joy2} we
find the equations are satisfied identically if $f_j$ is a holomorphic function with
respect to the complex structure on $M$. So we can construct a hyper-paracomplex
structure, with the opposite orientation, out of four holomorphic functions $f_1,...,f_4$.

Taking $f_1=f, f_2=f_3=0, f_4=1$ we obtain a local hyper-paracomplex structure. Consider
a particular neutral metric $g\in[g]$ such that

\centerline{$g(V_1,V_1)=g(V_2,V_2)=
-g(V_3,V_3)=-g(V_4,V_4)=1, g(V_j,V_k)=0, j\not= k$.}
The corresponding common Lee form is
given by

\centerline{$\theta = \frac{1}{f}\frac{\partial f}{\partial z^2}dz^2 +
 \frac{1}{\bar f}\frac{\bar{\partial} \bar f}{\bar{\partial} \bar z^2}dz^2$.}
Then $d\theta\not=0$ provided $\frac{\partial f}{\partial z^1}\not=0.$

\end{example}

\subsection{Twistor and reflector spaces on paraquaternionic K\"ahler manifold}
Consider the space $\tilde H_1$ of imaginary para-quaternions. It
is isomorphic to the Lorentz space $\mathbb R^2_1$ with a Lorentz
metric of signature (+,+,-) defined by
$<q,q'>=-Re(q\overline{q'})$, where $\overline q=-q$ is the
conjugate imaginary paraquaternion. In $\mathbb R^2_1$ there are
two kinds of 'unit spheres', namely the pseudo-sphere $S^2_1(1)$
of radius $1$ (the 1-sheeted hyperboloid) which consists of all
imaginary para-quaternions of norm $1$ and the pseudo-sphere
$S^2_1(-1)$ of radius (-1) (the 2-sheeted hyperboloid) which
contains all imaginary para-quaternions of norm (-1). The
1-sheeted hyperboloid $S^2_1(1)$ carries a natural paraHermitian
structure while the 2-sheeted hyperboloid $S^2_1(-1)$ carries a
natural Hermitian structure of signature (1,1), both induced by
the restriction of the Lorentz metric and the cross-product on
$\tilde H_1\cong \mathbb R^2_1$ defined by $$X \times
Y=\sum_{i\not= k}x^iy^kJ_iJ_k $$ for vectors $X=x^iJ_i,\quad
Y=y^kJ_k$. Namely, for a tangent vector $X=x^iJ_i$ to the
1-sheeted hyperboloid $S^2_1(1)$ at a point $q_+=q_+^kJ_k$ (resp.
tangent vector $Y=y_-^kJ_k$ to the 2-sheeted hyperboloid
$S^2_1(-1)$ at a point $q_-=q_-^kJ_k$) we define $PX:=q_+\times X$
(resp. $JY=q_-\times Y$). It is easy to check that $PX$ is again
tangent vector to $S^2_1(1), \quad P^2X=X, \quad <PX,PX>=-<X,X>$
(resp. $JY$ is tangent vector to $S^2_1(-1), \quad J^2Y=-Y, \quad
<JY,JY,>=<Y,Y,>$).

We start with a paraquaternionic K\"ahler manifold  $(M,g,\tilde
H=(J_{\alpha}))$ . The vector bundle $P$ carries a natural Lorentz
structure of signature (+,+,-) such that $(J_1,J_2,J_3)$ forms an
orthonormal local basis of $P$. The are two kinds of "unit sphere"
bundles according to the existence of the 1-sheeted hyperboloid
$S^2_1(1)$ and the 2-sheeted hyperboloid $S^2_1(-1)$. The twistor
space $Z^+(M)$ (resp. $Z^-(M)$) is the unit pseudo-sphere bundle
with fibre $S^2_1(1)$ (resp. $S^2_1(-1)$). In other words, the
fibre of $Z^+(M)$ consists of all almost paracomplex structures
(resp. all almost complex structures) compatible with the given
paraquaternionic K\"ahler structure. The bundle $Z^+(M)$ over a
4-dimensional manifold with a neutral metric was constructed in
\cite{JR} and called there \emph{the reflector space}. Further, we
keep their notation.

Denote by $\pi^{\pm}$ the projectionof $Z^{\pm}(M)$ onto $M$,
respectively. Keeping in mind the formal similarity with the
quaternionic geometry where there are two natural almost complex
structures \cite{AHS,ES}, we observe the existence of two
naturally arising  almost paracomplex structures on $Z^+(M)$
(resp. two almost complex structures on $Z^-(M)$ \cite{BDM})
defined as follows:

The Levi-Civita connection on $P$ preserves the Lorentz metric and
induces a linear connection on $Z^{\pm}$ i.e a splitting of the
tangent bundle $TZ^{\pm}=\mathbb H^{\pm}\otimes \mathbb V^{\pm}$,
respectively, where $\mathbb V^{\pm}$ is the vertical distribution
tangent to the fibre $S^2_1(1)$, (resp. $S^2_1(-1)$) and $\mathbb
H^{\pm}$ a supplementary horizontal distribution induced by the
Levi-Civita connection. By definition, the horizontal transport
associated to $\mathbb H^{\pm}$ preserves the canonical Lorentz
metric of the fibres $S^2_1(1)$ (resp. $S^2_1(-1)$); and also
their orientation; as a corollary, it preserves the canonical
paracomplex structure on $S^2_1(1)$ (resp. the canonical complex
structure on $S^(2_1(-1)$) described above. Since the vertical
distribution $\mathbb V^{\pm}$ is tangent to the fibres, this
paracomplex structure (resp. complex structure) induces an
endomorphism $\tilde P$ with $\tilde P^2=id$ (resp. $\tilde J$
with $\tilde J^2=-id$)  on $\mathbb V_z^{\pm}$ for each $z\in
Z^+(M)$ (resp. $Z^-(M)$). On the other hand, each point $z$ on
$Z^+(M)$ (resp. $Z^-(M)$) is by definition a paracomplex structure
on $T_{\pi(z)}Z^+(M)$ (resp. a complex structure on
$T_{\pi(z)}Z^-(M)$) which may be lifted into an endomorphism $\bar
P$ on $\mathbb H^+_z$ with $\bar P^2=id$ (resp. $\bar J$ on
$\mathbb H^-_z$ with $\bar J^2=-id$). We define an almost
paracomplex structures $\mathbb P_1,\mathbb P_2$ on $Z^+(M)$ and
an almost complex structure $\mathbb J_1,\mathbb J_2$ on $Z^-(M)$
by
\begin{gather*}
\mathbb P_1(\mathbb V^+)=\mathbb V^+, \quad \mathbb
P_1\vert_{\mathbb V^+}=\tilde P,\qquad \mathbb P_1(\mathbb
H^+)=\mathbb H^+, \quad \mathbb P_1\vert_{\mathbb H^+}=\bar P,\\
\mathbb P_2(\mathbb V^+)=\mathbb V^+, \quad \mathbb
P_2\vert_{\mathbb V^+}=-\tilde P,\qquad \mathbb P_2(\mathbb
H^+)=\mathbb H^+, \quad \mathbb P_2\vert_{\mathbb H^+}=\bar P;\\
\mathbb J_1(\mathbb V^-)=\mathbb V^-, \quad \mathbb
J_1\vert_{\mathbb V^-}=\tilde J,\qquad \mathbb J_1(\mathbb
H^-)=\mathbb H^-, \quad \mathbb J_1\vert_{\mathbb H^-}=\bar J,\\
\mathbb J_2(\mathbb V^-)=\mathbb V^-, \quad \mathbb
J_2\vert_{\mathbb V^-}=-\tilde J,\qquad \mathbb J_2(\mathbb
H^-)=\mathbb H^-, \quad \mathbb J_2\vert_{\mathbb H^-}=\bar J.
\end{gather*}

Define a pseudo Riemannian  metrics on $Z^+(M)$ (resp. $Z^-(M)$)
by $h_t^+=\pi^*g+t<,>^+_v, t\not =0, <,>^+_v$ being the
restriction of the Lorentz metric to the fibres $S^2_1(1)$ (resp.
$h_t^-=\pi^*g+t<,>^-_v, t\not =0, <,>^-_v$ being the restriction
of the Lorentz metric to the fibres $S^2_1(-1)$). It is easy to
check that $h^+$ (resp. $h^-$) is compatible with both $\mathbb
P_1, \mathbb P_2$ (resp. $\mathbb J_1, \mathbb J_2$) such that
$(Z^+(M),h^+_t, \mathbb P_{1,2})$ become  an almost paraHermitian
manifolds (resp. $(Z^-(M),h^-_t, \mathbb J_{1,2})$ become an
almost Hermitian manifolds).

The almost paracomplex structures $\mathbb P_1, \mathbb P_2$ and
the neutral metrics $h_t^+$ on the reflector space of a
4-dimensional manifold with a neutral metric $g$ are investigated
in \cite{JR}. The authors show that the almost paracomplex
structure $\mathbb P_2$ is never integrable while the almost
paracomplex structure $\mathbb P_1$ is integrable if and only if
the neutral metric $g$ is self dual. They also prove that the
neutral metric $h_t^+$ on the reflector space is Einstein if and
only if $g$ is self-dual Einstein and either $ts=12$ or $ts=6$.

Almost Hermitian geometry of $(Z^-(M),h^-_t, \mathbb J_{1,2})$ is
investigated in \cite{BDM}. The calculations there are completely
applicable to the almost paraHermitian geometry of $(Z^+(M),h^+_t,
\mathbb P_{1,2})$. In terms of the almost paraHermitian geometry
of $(Z^+(M),h^+_t, \mathbb P_{1,2})$ Theorem~1 and Theorem~2 in
\cite{BDM} read as follows
\begin{thm}\label{bdo1}
On the reflector space $(Z^+(M))$ of a paraquaternionic K\"ahler manifold
of dimension $4n\ge 8$ we have:
\begin{enumerate}
\item[i)] The almost paracomplex structure $\mathbb P_1$ is integrable and the Lee form of
the paraHermitian structure $(\mathbb P_1,h_t^+)$ is zero. The structure
$(\mathbb P_1,h_t^+)$ is paraK\"ahler if and only if $ts=4n(n+2)$;
\item[ii)] The almost paracomplex structure $\mathbb P_2$ is never integrable and
the Lee form of the almost paraHermitian structure $(\mathbb
P_2,h_t^+)$ is zero. The structure $(\mathbb P_2,h_t^+)$ is nearly
paraK\"ahler if and only if  $ts=2n(n+2)$ and  almost paraK\"ahler
if and only if $ts=-4n(n+2)$.
\end{enumerate}
\end{thm}
\begin{thm}\label{bdo2}
On the reflector space $(Z^+(M))$ of an oriented 4-dimensional manifold
 $M$ with a neutral metric $g$ we have the following:
\begin{enumerate}
\item[i)] The almost paraHermitian structure $(\mathbb P_1,h_t^+)$ has zero Lee form if and only if
the metric $g$ is self-dual.  It is paraK\"ahler if and only if the
metric $g$ is Einstein self-dual and $ts=12$;
\item[ii)] The Lee form of
the almost paraHermitian structure $(\mathbb P_2,h_t^+)$ is zero.
The structure $(\mathbb P_2,h_t^+)$ is  nearly paraK\"ahler if and
only if the metric $g$ is self-dual Einstein and $ts=6$ and almost
paraK\"ahler if and only if $ts=-12$.
\end{enumerate}
\end{thm}

\begin{rem}
On a  paraquaternionic manifold of dimension $4n\ge 8$ we may
construct the almost paracomplex structure $\mathbb P_1$ on the
reflector space $Z^+$ and the almost complex structure $\mathbb
J_1$ on the twistor space $Z^-$ using the horizontal distribution
generated by a torsion-free  connection instead of the horizontal
distribution of the Levi-Civita connection. In that case, we find
an  analogue of the result of S.Salamon \cite{S1,S2,S3}, (proved
also independently by L. Berard-Bergery, unpublished, see
\cite{Bes}). Namely, we have
\begin{thm}
On a paraquaternionic  manifold of dimension $4n\ge 8$
the almost paracomplex structure $\mathbb P_1$ on $Z^+$ and the almost complex
structure $\mathbb J_1$ on $Z^-$ are always integrable
\end{thm}
We sketch a proof which is completely similar to the proof in the
case of quaternionic  manifold presented in \cite{Bes}. Denote by
$R$  the curvature of a torsion-free  connection $\nabla$. Let
$S=xJ_1+yJ_2+zJ_3$ be either an almost paracomplex structure or an
almost complex structure compatible with the given
paraquaternionic structure i.e. the triple $(x,y,z)$ satisfies
either $x^2+y^2-z^2=1$ or $x^2+y^2-z^2=-1$. Denote by $$
S(R)(X,Y)=[R(SX,SY),S] - S[R(SX,Y),S] -
S[R(X,SY),S]+S^2[R(X,Y),S]. $$ In view of the analogy with the
proof in the quaternionic case presented in \cite{Bes},
14.72-14.74 the result will follow if $S(R)=0$. The last identity
can be checked in the exactly same way as it is done in
\cite{Bes}, Lemma~14.74 using \eqref{rel1} instead of formulas
14.39 in \cite{Bes}.
\end{rem}
\subsection{Examples of nearly paraK\"ahler and almost paraK\"ahler manifolds}
Theorem~\ref{bdo1} helps to find examples of nearly paraK\"ahler
and almost paraK\"ahler manifolds. We note that the sign of the
scalar curvature (if not zero) is not a restriction
 since the metric (-g) have
scalar curvature with opposite sign. Hence, taking the reflector space of any neutral (anti)
self-dual Einstein manifold with non-zero scalar curvature in dimension four and
any quaternionic paraK\"ahler manifold in dimension $4n\ge 8$
 we can find a real number $t$ to get nearly paraK\"ahler and almost paraK\"ahler
 structure on it.
\begin{enumerate}
\item The pseudo-sphere $S^3_3$ is endowed with an almost paracomplex structure
\cite{Lib} and it is shown in \cite{Bej1} that there exists a nearly paraK\"ahler
structure on $S^3_3$ induced from the so-called \emph{second kind Cayley numbers} (see
\cite{Lib}) in $\mathbb R^4_3$. The structure is Einstein with non-zero scalar curvature,
in fact the metric is the standard neutral metric on $S^3_3$ inherited from $\mathbb R^4_3$.
\item Start with one of the following 4-dimensional neutral self-dual Einstein spaces \\
$(S^2_2=SO^+(2,3)/GL^+(2,\mathbb R),can)$, $(\mathbb
CP^{1,1}=(SU(2,1)/(SO(1,1).U(1)),can)$, or $(SL(3,\mathbb
R)/GL^+(2,\mathbb R), c.Kill|_{sl(3,\mathbb R}))$, where $c$ is a
suitable constant and $Kill|_{sl(3,\mathbb R)}$ is the restriction
of the Killing form of $sl(3,\mathbb R)$  to the homogeneous space
\\ $SL(3,\mathbb R)/GL^+(2,\mathbb R$.  The corresponding
reflector spaces are $SO^+(2,3)/GL^+(2,\mathbb R)$,
$SU(2,1)/(SO(1,1).U(1))$, $SL(3,\mathbb R)/(\mathbb R^+\times
\mathbb R^+\cup \mathbb R^-\times\mathbb R^-)$, respectively
\cite{Kath}. These homogeneous spaces admit a homogeneous nearly
paraK\"ahler structure of non-zero scalar curvature as well as an
homogeneous  almost paraK\"ahler structure according to
Theorem~\ref{bdo2}
\item Non-homogeneous example arises from the non-(locally) homogeneous neutral self-dual
Einstein space of non-zero scalar curvature described in
\cite{BVq}. Its reflector space admit a nearly paraK\"ahler
structure of non-zero scalar curvature, as well as an almost
paraK\"ahler structure, due to Theorem~\ref{bdo2}
\end{enumerate}

To the best of our knowledge there are no known examples of Ricci
flat 6-dimensional Nearly paraK\"ahler manifolds.

\bibliographystyle{hamsplain}

\begin{thebibliography}{12}

\bibitem{Akiv} M. Akivis, {\em On the real theory of conformal structures},
J. Geom. Phys. {\bf 21} (1996), 55-80.

\bibitem{An} A. Andrada, {\em Complex product structures on differentiable manifolds},
preprint.

\bibitem{Andr} A. Andrada, {\em Hypersymplectic four-dimensional Lie algebras},
math.DG/0310462.

\bibitem{ASal} A. Andrada, S. Salamon, {\em Complex product structures on Lie algebras},
Forum Math. {\bf 17} (2005), no. 2, 261--295. 

\bibitem{Apos1} V. Apostolov, {\em Generalized Goldberg-Sachs theorems
for pseudo-Riemannian four-manifolds}, J. Geom. Phys. {\bf 27} (1998), 185-198.

\bibitem{AG1} V. Apostolov, P. Gauduchon, {\em The Riemannian Goldberg-Sachs
theorem}, Int. J. Math. {\bf 8} (1997), 421-439.

\bibitem{Ash} A. Ashtekar, T. Jacobson, L. Smolin, {\em A new characterization of half-flat
solutions to Einstein's equation}, Comm. Math. Phys. {\bf 115} (1988), 631-648.

\bibitem{AHS} M.F. Atiyah, N.J. Hitchin, I.M. Singer, {\em Self-duality
in four-dimensional Riemannian geometry}, Proc. Roy. Soc. London Ser. A
{\bf 362} (1978), 425-461.

\bibitem{BE} T. Bailey, M. Eastwood, {\em Complex paraconformal manifolds-their
differential geometry and twistor theory},  Forum Math. {\bf 3} (1991), 61-103.

\bibitem{Bal} V. Balashchenko, {\em Invariant nearly Kähler $f$-structures on
homogeneous spaces},
 Global differential geometry: the mathematical legacy of Alfred Gray (Bilbao,
   2000), 263--267, Contemp. Math., {\bf 288}, Amer. Math. Soc., Providence, RI, 2001.

\bibitem{Bar} J. Barret, G.W. Gibbons,M.J. Perry, C.N. Pope, P.Ruback,
{\em Kleinian geometry and the N=2 superstring}, Int. J. Mod. Phys. {\bf A9}
(1994), 1457-1494.

\bibitem{Bej} C. Bejan, {\em A classification of the almost paraHermitian
manifolds}, Proc. Conf. 'Diff. Geom. Appl.' (Dubrownik, 1988), 23-27,
Univ. Novi Sad, 1989.

\bibitem{Bej1}  C. Bejan, {\em Some examples of manifolds with hyperbolic
structures}, Rend. Math., Ser. VII {\bf 14}, Roma (1994), 557-565.

\bibitem{BO} C. Bejan, L. Ornea, {\em An example of an almost hyperbolic
Hermitian manifold}, Internat. J. Math. Math. Sci. {\bf 21}, (1998), 613--618.

\bibitem{BM} F. Belgun, A. Moroianu, {\em Nearly K\"ahler manifolds with
reduced holonomy}, Ann. Glob. Anal. Geom. {\bf 19} (2001), 307-319.

\bibitem{Bes} A. Besse, Einstein manifolds, Ergebnisse der Mathematik und
ihrer Grenzgebiete 3. Folge, Band 10, Springer, Berlin, 1987.

\bibitem{Biel} R. Bielawski, {\em Manifolds with an $SU(2)$-action on the
tangent bundle}, math.DG/0309301.

\bibitem{Biq} O. Biquard, {\em Les équations de Seiberg-Witten sur une
surface complexe non kählérienne},
Comm. Anal. Geom. {\bf 6} (1998), 1, 173--197.

\bibitem{BDM} D. Blair, J. Davidov, O. Muskarov, {\em Hyperbolic twistor spaces},
Rocky Mountain J. Math., in press.

\bibitem{BL1} N. Blazic, {\em  Projective space and
pseudo-Riemannian geometry}, Publ. Inst. Math. {\bf 60} (1996), 101-107.

\bibitem{BVq} N. Blazic, S. Vukmirovic, {\em Examples of self-dual, Einstein
metrics of (2,2)-signature},  Math. Scand.  {\bf 94}  (2004),  no.
1, 63--74.

\bibitem{BVuk0}  N. Blazic, S. Vukmirovic, {\em Para-hypercomplex structures on a four-dimensional Lie group}, 
Contemporary geometry and related topics, 41--56, World Sci. Publishing, River Edge, NJ, 2004. 

\bibitem {Bo} C. Boyer, {\em A note on hyperhermitian four manifolds},
Proc. Amer. Math. Soc. ${\bf 102}$ No 1 $(1988)$, $157-164$.

\bibitem{CFG} V. Cruceanu, P. Fortuny, P.M. Gadea, {\em
A survey on paracomplex geometry}, Rocky Mountain J. Math.
{\bf 26} (1996), 83--115.

\bibitem{DS}V. De Smedt, S. Salamon, {\em Anti-self-dual metrics on Lie groups},
Proc. Conf. Integrable Systems and Differential Geometry, Contemp. Math. {\bf 308}
(2002), 63-75.

\bibitem{D1} M. Dunajski, {\em Hyper-complex four-manifolds from Tzitzeica
equation}, J. Math. Phys. {\bf 43} (2002), 651-658.

\bibitem{ES} J. Eells, S.Salamon, {\em Twistorial construction of harmonic maps
of surfaces into four-manifolds}, Ann. Scuola Norm. Sup. Pisa Cl. Sci.
{\bf 12} (1985), 589-640.

\bibitem{erd} S. Erdem, {\em Paraholomorphic structures and the connections
of vector bundles over paracomplex manifolds} New Zealand J. Math.
{\bf 30} (2001), 41-50.

\bibitem{FPed} A. Fino, H. Pedersen, Y.-S. Poon, M.W. Sorensen,
{\em Neutral Calabi-Yau structures on Kodaira manifolds},  Comm.
Math. Phys.  {\bf 248}  (2004),  no. 2, 255--268.

\bibitem{GM} P.M. Gadea, J.M. Masque, {\em Classification of the almost
paraHermitian manifolds}, Rend. Mat. Roma {\bf 11} (1991), 377-396.

\bibitem{GRio} E. Garcia-Rio,Y. Matsushita, R. Vasquez-Lorentzo, {\em
Paraquaternionic K\"hler manifold}, Rocky Mountain J. Math. {\bf 31}
(2001), 237-260.

\bibitem{Gau} P. Gauduchon, {\em La $1$-forme de torsion d'une variété
hermitienne compacte}, Math. Ann. {\bf 267} (1984),  495--518.

\bibitem{Ga1}  P. Gauduchon, {\em Hermitian connections and Dirac operators},
 Boll. Un. Mat. Ital. {\bf B (7) 11} (1997), no. 2, suppl., 257--288.

\bibitem{GT} P. Gauduchon, P.K.Tod, {\em Hyper-Hermitian metrics with symmetry},
J. Geom. Phys. {\bf 25} (1998), 291-304.

\bibitem{Gr1} A. Gray, {\em Nearly K\"ahler manifolds}, J. Diff. Geom
{\bf 4} (1970), 283-309.

\bibitem{Gr2} A. Gray, {\em Riemannian manifolds with geodesic symmetries of
order 3}, J. Diff. Geom {\bf 6} (1972), 343-369.

\bibitem{Gr3} A. Gray, {\em The structure of Nearly K\"ahler manifolds},
Math. Ann. {\bf 223} (1976), 233-248.

\bibitem{KH} K. Hasegawa, {\em Four-dimensional compact
solvmanifolds with and without compex analytic structures},
math.CV/0401413.

\bibitem{Hit} N. Hitchin, {\em Hypersymplectic quotients}, Acta Acad. Sci.
Tauriensis {124} supl., (1990), 169-180.

\bibitem{hul} C.M. Hull, {\em Actions for (2,1)sigma models and strings},
Nucl. Phys. {\bf B 509} (1988), no.1, 252-272.

\bibitem{Hull} C.M. Hull, {\em A Geometry for Non-Geometric String
Backgrounds}, hep-th/0406102.

\bibitem{IT} S. Ivanov, V. Tsanov, {\em Complex product structures on some simple Lie
groups}, math.DG/0405584.

\bibitem{ITZ} S. Ivanov, V. Tsanov, S. Zamkovoy,
{\em Hyper-ParaHermitian manifolds with torsion}, J. Geom. Phys., in press, math.DG/0405585.

\bibitem{JR} G.Jensens, M. Rigoli, {\em Neutral surfaces in neutral four spaces},
Mathematische (Catania) {\bf 45} (1990), 407-443.

\bibitem{Joy} D. Joyce, {\em Explicit construction of self-dual 4-manifolds},
Duke Math. J. {\bf 77} (1995), 519-552.

\bibitem{Kam} H. Kamada, {\em Neutral hyper-K\"ahler structures on primary
Kodaira surfaces}, Tsukuba J. Math. {\bf 23} (1999), 321-332.

\bibitem{KK} S. Kaneyuki, M. Kozai, {\em Paracomplex structures and affine
symmetric spaces}, Tokyo J. Math. {\bf 8} (1985), 81-98.

\bibitem{Kath} I. Kath, {\em Killing spinors on pseudo-Riemannian Manifolds},
Habilitationscgrift Humboldt Universit\"at Berlin, 1999.

\bibitem {Ka} M. Kato, {\em Compact Differentiable 4-folds with Quaternionic
Structures}, Math. Ann. ${\bf 248}$ $(1980)$, $79-96$.

\bibitem{Kir} V. Kirichenko, {\em K-spaces of maximal rank}, (in russian),
Mat. Zam. {\bf 22} (1977), 465-476.

\bibitem{Lib} P. Libermann, {\em Sur le probleme d'equivalence de certains
structures infinitesimales}, Ann. Mat. Pura Appl. {\bf 36} (1954), 27-120.

\bibitem{Mal} A. I. Malcev, {\em On a class of homogeneous spaces},
reprinted in Amer. Math. Soc. Translations, Series 1, {\bf 9} (1962), 276-307.

\bibitem{Ma} Y. Matsushita, {\em Fields of $2$-planes and two kinds of almost
complex structures on compact $4$-dimensional manifolds} Math. Z. {\bf 207}
 (1991), 281-291.

\bibitem{Mil} J. Milnor, {\em Curvatures of left invariant metrics on Lie groups},
 Advances in Math. {\bf 21} (1976), no. 3, 293-329.

\bibitem{OV} H. Ooguri, C. Vafa, {\em Geometry of $N=2$ strings},
Nucl. Phys. {\bf B 361} (1991), 469-518.

\bibitem{Ol} Z. Olszak, {\em Four-dimensional paraHermitian manifold},
Tensor N.S. {\bf 56} (1995), 215-226.

\bibitem{PR} R.Penrose, W. Rindler, Spinors and Space time 2, CUP, 1986.

\bibitem{Pet} J. Petean, {\em Indefinite K\"ahler metrics on compact
complex surfaces}, Comm. Math. Phys. {\bf 189} (1997), 227-235.

\bibitem{RV} F. Raymond, A. Vasquez, {\em
$3$-manifolds whose universal coverings are Lie groups},
Topology Appl. {\bf 12} (1981), 161-179.

\bibitem{S1} S. Salamon, {\em Quaternionic Kähler manifolds},
Invent. Math. {\bf 67} (1982), 143-171.

\bibitem{S2} S. Salamon, {\em Quaternionic manifolds},  Symposia Mathematica,
{\bf 26} (1982), 139-151.

\bibitem{S3} S. Salamon, {\em Differential geometry of quaternionic manifolds},
 Ann. Sci. École Norm. Sup. {\bf (4) 19} (1986), 31-55.

\bibitem{Scot} P. Scott, {\em The geometries of 3-manifolds}, Bull. London Math.
Soc. {\bf 15} (1983), 407-487.

\bibitem{Vuk} S. Vukmirovic, {\em Paraquaternionic  reduction}, math.DG/0304424.

\bibitem{Wall} C.T.C. Wall, {\em Geometric structures on compact complex analytic
surfaces}, Topology {\bf 25} (1986), 119-153.

\bibitem{Yano} K. Yano, {\em Affine connections in an almost product spaces},
Kodai Math. Sem. Rep. {\bf 11} (1959), 1-24.

\bibitem{Yano1} K. Yano, Differential geometry on complex and almost complex
spaces. International Series of Monographs in Pure and Applied Mathematics,
Vol. 49 A Pergamon Press Book. The Macmillan Co., New York 1965.


\end{thebibliography}

\providecommand{\bysame}{\leavevmode\hbox to3em{\hrulefill}\thinspace}






\end{document}